\documentclass[12pt]{amsart}
\font\emailfont=cmtt10

\usepackage{amsmath,amsthm,amsfonts,amscd,flafter,epsf}

\newcommand\commentable[1]{#1}

\newcommand{\Tors}{\mathrm{Tors}}

\newcommand{\HF}{HF}

\newtheorem{theorem}{Theorem}[section]
\newtheorem{prop}[theorem]{Proposition}
\newtheorem{cor}[theorem]{Corollary}

\newtheorem{lemma}[theorem]{Lemma}

\newtheorem{defn}[theorem]{Definition}

\def\endproof{\relax\ifmmode\expandafter\endproofmath\else
  \unskip\nobreak\hfil\penalty50\hskip.75em\hbox{}\nobreak\hfil\bull
  {\parfillskip=0pt \finalhyphendemerits=0 \bigbreak}\fi}
\def\endproofmath$${\eqno\bull$$\bigbreak}
\def\bull{\vbox{\hrule\hbox{\vrule\kern3pt\vbox{\kern6pt}\kern3pt\vrule}\hrule}}

\newcommand{\Q}{\mathbb{Q}}
\newcommand{\R}{\mathbb{R}}

\newcommand{\Z}{\mathbb{Z}}

\newcommand{\Zmod}[1]{\Z/{#1}\Z}

\newcommand{\Ker}{\mathrm{Ker}}

\newcommand{\cm}{\cdot}

\newcommand{\Nbd}[1]{{\mathrm{nd}}(#1)}
\newcommand{\nbd}[1]{\Nbd{#1}}

\newcommand{\ModSWfour}{\mathcal{M}}
\newcommand{\ModFlow}{\ModSWfour}

\newcommand{\SpinC}{{\mathrm{Spin}}^c}

\newcommand\Wedge{\Lambda}

\newcommand\abuts\Rightarrow
\newcommand\Sym{\mathrm{Sym}}

\newcommand\spinccan{\ell}

\newcommand\HFpRed{\HFp_{\red}}

\newcommand\sRelSpinC{\underline\spinc}

\newcommand\Filt{\mathcal F}

\newcommand\x{\mathbf x}
\newcommand\w{\mathbf w}

\newcommand\y{\mathbf y}

\newcommand\ModSphere{\ModFlow\left({\mathbb S}\longrightarrow 
\Sym^{g-1}(\Sigma_{1})\times \Sym^2(\Sigma_{2})\right)}
\newcommand\ModSpheres\ModSphere
\newcommand\CF{CF}

\newcommand\CFa{\widehat{CF}}
\newcommand\CFp{\CFb}
\newcommand\CFm{\CF^-}

\newcommand\HFleq{\HF^{\leq 0}}

\newcommand{\red}{\mathrm{red}}

\newcommand\HFp{\HFb}

\newcommand\HFm{\HF^-}
\newcommand\CFinf{CF^\infty}
\newcommand\HFinf{HF^\infty}
\newcommand\CFb{CF^+}
\newcommand\HFa{\widehat{HF}}
\newcommand\HFb{HF^+}
\newcommand\gr{\mathrm{gr}}
\newcommand\Mas{\mu}
\newcommand\UnparModSp{\widehat \ModSp}
\newcommand\UnparModFlow\UnparModSp
\newcommand\Mod\ModSp

\newcommand{\cald}{{\mathcal D}}

\newcommand\PD{\mathrm{PD}}

\newcommand{\spinc}{\mathfrak s}

\newcommand{\spinct}{\mathfrak t}

\newcommand\ModMaps{\mathcal M}
\newcommand\ModSp\ModMaps

\newcommand\Ta{{\mathbb T}_{\alpha}}
\newcommand\Tb{{\mathbb T}_{\beta}}
\newcommand\Tc{{\mathbb T}_{\gamma}}
\newcommand\Td{{\mathbb T}_{\delta}}

\newcommand\alphas{\mbox{\boldmath$\alpha$}}

\newcommand\betas{\mbox{\boldmath$\beta$}}
\newcommand\gammas{\mbox{\boldmath$\gamma$}}
\newcommand\deltas{\mbox{\boldmath$\delta$}}

\newcommand\uHFp{\underline{\HF}^+}

\newcommand\uCFp{\underline{\CFp}}

\newcommand\Fp[1]{F^{+}_{#1}}

\newcommand\spincu{\mathfrak u}

\newcommand\Field{\mathbb F}

\newcommand\HFc{\HF^\circ}
\newcommand\CFc{\CF^\circ}

\newcommand\Dual{\mathcal D}
\newcommand\Duality\Dual

\newcommand\InjMod[1]{{\mathcal T}^+_{#1}}

\newcommand\MCone{M}

\newcommand\CFK{CFK}

\newcommand\CFKinf{\CFK^{\infty}}

\newcommand\rPhi{\overline \phi}
\newcommand\rC{\overline C}
\newcommand\Width{b}
\newcommand\imult{\mathfrak m}

\newcommand\Xa{\widehat{\mathbb X}}
\newcommand\Xp{{\mathbb X}^+}
\newcommand\Xd{{\mathbb X}^{\delta}}

\newcommand\vertYa{\widehat{\mathfrak v}}
\newcommand\horYa{\widehat{\mathfrak h}}
\newcommand\AYa{\widehat{\mathfrak A}}
\newcommand\BYa{\widehat{\mathfrak B}}
\newcommand\Aa{\widehat{A}}
\newcommand\Ba{\widehat{B}}

\newcommand\verta{\widehat v}
\newcommand\hora{\widehat h}
\newcommand\vertp{v^+}
\newcommand\horp{h^+}
\newcommand\vertd{v^{\delta}}
\newcommand\hord{h^{\delta}}

\newcommand\fp{f^+}
\newcommand\Hp[1]{H^+_{#1}}

\newcommand\Ap{{A}^+}
\newcommand\Bp{{B}^+}

\newcommand\Ad{{A}^{\delta}}
\newcommand\Bd{{B}^{\delta}}

\newcommand\fa{\widehat f}

\newcommand\fd{f^{\delta}}
\newcommand\Hd{H^{\delta}}

\newcommand\CFd{CF^{\delta}}
\newcommand\HFd{HF^{\delta}}
\newcommand\Dd{D^{\delta}}

\newcommand\BigAp{{\mathbb A}^+}
\newcommand\BigBp{{\mathbb B}^+}
\newcommand\Dp{D^+}

\newcommand\CapSurf{\widehat F}

\newcommand\BigAd{{\mathbb A}^{\delta}}
\newcommand\BigBd{{\mathbb B}^{\delta}}
\newcommand\spincx{\mathfrak x}
\newcommand\spincy{\mathfrak y}

\newcommand\BigAa{\widehat{\mathbb A}}
\newcommand\BigBa{\widehat{\mathbb B}}

\commentable{

\title[{Knot Floer homology and integer surgeries}] 
{Knot Floer homology and integer surgeries}

\author[Peter Ozsv{\'a}th]{Peter Ozsv\'ath}
\address{Department of
Mathematics, Columbia University, New York 10027 \newline
\indent{\emailfont{petero@math.columbia.edu}}}
\thanks{PSO was supported by NSF grant number DMS 0234311}

\author[Zolt{\'a}n Szab{\'o}]{Zolt{\'a}n Szab{\'o}}
\address{Department of Mathematics, Princeton University, New Jersey
  08544 \newline \indent{\emailfont{szabo@math.princeton.edu}}}}
\thanks{ZSz was supported by NSF grant number DMS 0107792}

\begin{document}

\begin{abstract}  
  Let $Y$ be a closed three-manifold with trivial first homology, and
  let $K\subset Y$ be a knot.  We give a description of the Heegaard
  Floer homology of integer surgeries on $Y$ along $K$ in terms of the
  filtered homotopy type of the knot invariant for $K$.  As an
  illustration of these techniques, we calculate the Heegaard Floer
  homology groups of non-trivial circle bundles over Riemann surfaces
  (with coefficients in $\Zmod{2}$).
\end{abstract} 

\maketitle

\section{Introduction}

Heegaard Floer homology is an invariant for closed-oriented
three-manifolds $Y$~\cite{HolDisk}. The invariant, denoted
$\HFc(Y)$\footnote{In~\cite{HolDisk}, we defined several variants of
  Heegaard Floer homology -- $\HFa$, $\HFm$, $\HFinf$, and $\HFp$,
  which are related to one another by various exact sequences. We
  denote this entire collection here by $\HFc$. In this paper,
  however, we will focus primarily on the case of $\HFp$.}, is the
homology of a chain complex whose generators have a combinatorial
definition, and whose boundary operator counts certain
pseudo-holomorphic disks in associated spaces.  Moreover, if $W\colon
Y_1\longrightarrow Y_2$ is a smooth, connected, oriented cobordism
from the connected, oriented three-manifold $Y_1$ to $Y_2$, equipped
with a $\SpinC$ structure, there is an associated map of
Heegaard Floer homology groups.  Explicit calculations of the chain
complex and associated chain maps are difficult to make in general,
though under favorable circumstances, the homology groups can be
determined with the help of a long exact sequence relating the
Heegaard Floer homologies of various surgeries on a given knot,
cf.~\cite{HolDiskTwo}.

In~\cite{Knots} and~\cite{RasmussenThesis}, a closely related
invariant is defined for null-homologous knots $K$ in a closed,
oriented three-manifold $Y$, taking the form of an induced filtration
on the Heegaard Floer complex of $Y$. The filtered chain
homotopy type of this complex is a knot invariant, to which we
refer loosely as ``knot Floer homology''.

If $K\subset Y$ is null-homologous, there is a canonical
identification of framings on $K$ with integers. Given an integer $n$,
let $Y_n(K)$ denote the three-manifold obtained by $n$-framed surgery
on $Y$ along $K$.  When $n$ is sufficiently large, there is an
immediate relationship between the knot Floer homology of $K$ and the
Heegaard Floer homology of $Y_n(K)$, see~\cite{Knots},
\cite{RasmussenThesis}.  In cases where $Y$ is sufficiently simple --
for example, when $Y\cong S^3$ -- these data are sufficient to
determine $\HFp(Y_n(K))$ for arbitrary integers $n$. However, even
this case, one loses information on some of the additional structure
carried by Heegaard Floer homology (for example, its absolute
grading).

The aim of this paper is to describe completely the Heegaard Floer
homology of integral surgeries on $Y$ along $K$ in terms of data
associated to the knot $K$, in the case where $Y$ has trivial first
homology.  To describe this construction, we recall some aspects of
knot Floer homology.

Knot Floer homology associates to a knot $K$ in $Y$ a
$\Z\oplus\Z$-filtered $\Z[U]$-complex $C=\CFKinf(Y,K)$. More precisely,
$C$ is generated over $\Z$ by a set $X$ equipped with a function
$\Filt \colon X \longrightarrow \Z\oplus \Z$ with the property that if
$\x\in X$ is an element with $\Filt(\x)=(i,j)$, then $\partial \x$ can
be written as a linear combination of elements $\y\in X$ with
$\Filt(\y)\leq \Filt(\x)$.  Moreover, the action of $U$ on $C$ is
induced from an action of $U$ on the generating set $X$,
with the property that for each $\x\in X$ with $\Filt(\x)=(i,j)$, then
$\Filt(U\cm \x)=(i-1,j-1)$.

We write the two factors of the $\Z\oplus\Z$ filtration as
$(i,j)$. Let  $S$ be a region
in the plane which has the property that
for all $(i,j)\in S$ and all $(i',j')\geq (i,j)$, we have that 
$(i',j')\in S$. Then the subset of $C$ generated by points with filtration
level contained in $S$ naturally inherits the structure of a quotient
complex, which we write as $C\{S\}$. 
For example,  $C\{i\geq 0\}$ denotes the quotient complex
of $C$ generated by $\x\in X$ with $\Filt(\x)=(i,j)$ with $i\geq 0$.
There is a canonical (up to sign) chain homotopy equivalence $C\{i\geq
0\}\simeq C\{j\geq 0\}$. Indeed, $\Bp=C\{i\geq 0\}$ is identified with
$\CFp(Y)$.

Define $\Ap_s=C\{\max(i,j-s)\geq 0\}$. There are two canonical chain
maps $\vertp_s\colon \Ap_s \longrightarrow \Bp$ and $\horp_s\colon
\Ap_s\longrightarrow \Bp$. The map $\vertp_s$ is projection onto
$C\{i\geq 0\}$, while $\horp_s$ is projection onto $C\{j\geq s\}$,
followed by the identification with $C\{j\geq 0\}$ (induced by
multiplication by $U^s$), followed by the chain homotopy equivalence
from $C\{j\geq 0\}$ to $C\{i\geq 0\}$.

Let $\BigAp=\bigoplus_{s\in\Z} \Ap_s$ and $\BigBp=\bigoplus_{s\in\Z}
\Bp_s$ (for this latter group, each summand is isomorphic to $\Bp$,
but we include a subscript to distinguish the various summands), and
let $\Dp_n\colon \BigAp \longrightarrow \BigBp$ be the map
$$\Dp_n(\{a_s\}_{s\in\Z}) = \{b_s\}_{s\in\Z},$$ where here $$b_s =
\horp_{s-n}(a_{s-n})+\vertp_s(a_s).$$ Let $\Xp(n)$ denote the mapping
cone of $\Dp_n$; i.e. this is the chain complex whose underlying group
is $\BigAp\oplus\BigBp$, and whose differential over $\Zmod{2}$ has
the form $$\left(\begin{array}{rr}
\partial_{\BigAp} & 0 \\
\Dp_n & \partial_{\BigBp}
\end{array}\right).$$

The following result says that the above data associated to the knot
$K\subset Y$ contains enough information to deduce the Floer homologies of
all the three-manifolds obtained by integer surgeries on $K$, and all
the corresponding maps induced on homology by the natural cobordisms:

\begin{theorem}
\label{thm:IntegerSurgeries}
Let $Y$ be an integral homology three-sphere.  For any non-zero
integer $n$, the homology of the mapping cone $\Xp(n)$ of
$$\Dp_n\colon \BigAp\longrightarrow \BigBp$$
is isomorphic to
$\HFp(Y_n(K))$. Moreover, under this identification the natural map
$\HFp(Y)\cong H_*(\Bp_s)\longrightarrow H_*(\Xp(n))$ is identified with 
the map $\HFp(Y) \longrightarrow \HFp(Y_n(K))$ induced by the natural
two-handle cobordism from $Y$ to $Y_n(K)$ endowed with the $s^{th}$
$\SpinC$ structure (for some identification of these $\SpinC$
structures with $\Z$).
\end{theorem}

We have chosen our hypotheses of Theorem~\ref{thm:IntegerSurgeries} to
simplify the statement. For generalizations and refinements, see
below. In particular, see Theorem~\ref{thm:PreciseIntegerSurgeries}
for a more precise statement of the isomorphism, which takes gradings
into account, and Theorem~\ref{thm:Maps}, which makes
explicit the identification of the maps induced by cobordisms. See
also Subsection~\ref{subsec:ZeroSurgery} for a discussion of the case
where $n=0$, and Subsection~\ref{subsec:Generalize} for further
generalizations of Theorem~\ref{thm:IntegerSurgeries}.  We return to
the more general case of rational surgeries in the
sequel~\cite{KnotsRat}.

Note that knowing the filtered chain homotopy type of the knot
filtration uniquely determines the groups $\Ap_s$ and $\Bp_s$, and
also the maps $\vertp_s$. The maps $\horp_s$ are, in general, not
known explicitly, as they involve the chain homotopy equivalence
$C\{i\geq 0\}\simeq C\{j\geq 0\}$. Despite this shortcoming,
Theorem~\ref{thm:IntegerSurgeries} is still quite helpful in
performing explicit calculations, as we shall see in
Section~\ref{sec:Calculations}.  It is also worth noting that there is
an overall $\pm 1$ ambiguity in the homotopy equivalence, and hence in
the maps $\horp_s$; however, it is easy to see that the homology of
the resulting complex is independent of this ambiguity.

In Section~\ref{sec:Review}, we recall some of the aspects of Heegaard
Floer homology used later, and also set up some notation. In
Section~\ref{sec:ExactSeq}, we prove an exact sequence relating
Heegaard Floer homology groups of $Y_n(K)$, $Y_{m+n}(K)$, and $Y$,
which is an ingredient in the proof of
Theorem~\ref{thm:IntegerSurgeries}.  In Section~\ref{sec:Proof}, we
state and prove a more precise version of
Theorem~\ref{thm:IntegerSurgeries} (cf.
Theorems~\ref{thm:PreciseIntegerSurgeries} and~\ref{thm:Maps} below),
taking into account the gradings on Floer homology. In
Section~\ref{sec:Calculations}, we give some sample calculations to
illustrate the techniques from this paper.  As an example, we
calculate the reduced Heegaard Floer homology of any non-trivial
circle bundle over a Riemann surface (with coefficients in $\Zmod{2}$).

Theorem~\ref{thm:IntegerSurgeries}, of course, strengthens the
connection between Floer homology for knots and Floer homology for
closed three-manifolds. This could be further pursued from several
angles.  For example, this result could be viewed as motivation for
studying knot invariants in gauge-theoretic contexts, such as
Seiberg-Witten or Donaldson's theories (compare~\cite{CollinSteer}).
In a different direction, there seems to be a close connection between
knot Floer homology and Khovanov's homology for links
(cf.~\cite{Khovanov}, \cite{KhovanovRozansky}, \cite{BarNatan},
\cite{EunSooLee}), see also recent work of
Rasmussen~\cite{RasmussenSlice}. It is an open problem whether
Khovanov's homology admits an extension to a three-manifold invariant.

\vskip.2cm
\noindent{\bf{Acknowledgements.}}
The authors wish to thank Tomasz Mrowka, Jacob Rasmussen, and
Andr{\'a}s Stipsicz for valuable conversations during the course of
this work. We are especially grateful to Slaven Jabuka and Thomas Mark
for some very useful correspondence regarding the calculations
of circle bundles.

<\section{Review}
\label{sec:Review}

The purpose of this section is to introduce notation, recall some
elements of Heegaard Floer homology used later, and also to put in
place various preliminary notions.  In
Subsection~\ref{subsec:Gradings}, we set up terminology on gradings
which we will use throughout this paper; in
Subsection~\ref{subsec:ChainComplexes}, we recall some standard
terminology from homological algebra; in
Subsection~\ref{subsec:CFcx} we set up notation for the Heegaard
Floer complexes; in
Subsection~\ref{subsec:IntegralSurgeries}, we set up notation and
conventions for integral surgeries on null-homologous knots; in
Subsection~\ref{subsec:LargeNSurgeries} we recall the relationship
between knot Floer homology and Heegaard Floer homology of surgeries
with sufficiently large (integral) coefficients (compare~\cite{Knots}
and \cite{RasmussenThesis}); in Subsection~\ref{subsec:Unknot} we
verify Theorem~\ref{thm:IntegerSurgeries} for the unknot in $S^3$.
This calculation will be useful to us in the proof of
Theorem~\ref{thm:IntegerSurgeries}.  Finally, in
Subsection~\ref{subsec:Algebra}, we include some simple observations
about the $\Z[U]$-modules which we encounter in this paper.

\subsection{Gradings}
\label{subsec:Gradings}

Let $C$ be a free $\Z$ module which is freely generated by some set
$X$.  We say that $C$ is {\em relatively $\Z$-graded} if there is a
function $$\Delta\colon X \times X \longrightarrow \Z$$ with the
property that $$\Delta(x,y)+\Delta(y,z) = \Delta(x,z).$$ An element of
$C$ is said to be {\em homogeneous} if it can be written as a linear
combination of elements from a subset $S\subset X$ with the property
that for all $x,y\in S$, $\Delta(x,y)=0$.  An endomorphism $f\colon C
\longrightarrow C$ is said to be {\em homogeneous of degree $d$} if
the image of any homogeneous element of $C$ is homogeneous, and indeed
for all $x\in X$, $f(x)$ can be written as a linear combination of
elements $y$ with $\Delta(y,x) = d$. The same terminology can be used
when the relative grading takes values in the rational numbers, rather
than the integers (in which case we call it a {\em relative
$\Q$-grading}). A {\em relatively graded chain complex} is a
relatively graded group $C$, which is equipped with a differential
$$\partial \colon C \longrightarrow C$$ which is homogeneous of degree
$-1$.

Let $C$ be a free $\Z$-module which is freely generated by some set
$X$.  We say that $C$ is {\em absolutely graded} if $X$ is
equipped with a function $$\gr \colon X
\longrightarrow \Q.$$ An absolute grading $\gr$, of
course, induces a relative $\Q$-grading by the formula $$\Delta(x,y) =
\gr(x) - \gr(y).$$ In this case, the
absolute grading $\gr$ is said to be a lift of the
relative grading $\Delta$.
If $C_1$ and $C_2$ are graded chain complexes, then a map
$f\colon C_1 \longrightarrow C_2$ is said to be homogeneous of degree $c$
if it carries homogeneous elements in $C_1$ with degree $a$ 
to homogeneous elements in $C_2$ with degree $a+c$.

\subsection{Chain complexes}
\label{subsec:ChainComplexes}

We will often use standard notions from homological algebra, which we collect here for the reader's convenience.

Let $f\colon C_1 \longrightarrow C_2$ be a chain map between
$\Zmod{2}$-graded chain complexes.  Then, the {\em mapping cone}
$\MCone(f)$ is the chain complex whose underlying group is $C_1\oplus
C_2$, endowed with a differential 
$$(a,b) \mapsto (\partial a, \partial b + (-1)^{\gr(a)} \cm f(a)).$$

There is a short exact sequence of chain maps
$$
\begin{CD}
0@>>> C_2 @>{\iota}>> \MCone(f) @>{\pi}>> C_1 @>>> 0,
\end{CD}
$$
whose induced connecting homomorphism is identified (up to sign)
with the map on homology induced by $f$,
$$F_*\colon H_*(C_1) \longrightarrow H_*(C_2).$$

Two maps $f, f'\colon C_1\longrightarrow C_2$ are {\em chain
homotopic} if there is a map $h\colon C_1 \longrightarrow C_2$ with
$\partial\circ h - h\circ \partial = f-f'$. Two chain complexes $C_1$
and $C_2$ are said to be {\em chain homotopy equivalent} if there are
chain maps $\phi\colon C_1\longrightarrow C_2$ and $\psi\colon
C_2\longrightarrow C_1$ so that $\phi\circ \psi$ and $\psi\circ\phi$
are chain homotopic to the respective identity maps. The maps $\phi$
and $\psi$ are called {\em chain homotopy equivalences}. Chain homotopic
maps give rise to chain homotopy equivalent mapping cones.

Let $C_1$ and $C_2$ be a pair of chain complexes.  A {\em
quasi-isomorphism} is a chain map $f\colon C_1\longrightarrow C_2$
which induces an isomorphism in homology. Of course, a chain homotopy
equivalence is a quasi-isomorphism. Two complexes $C_1$ and $C_2$ are
said to be {\em quasi-isomorphic} if there is a third chain complex
$C_0$ and quasi-isomorphisms from $C_0$ to $C_1$ and $C_0$ to $C_2$.
The following lemma is standard:

\begin{lemma}
\label{lemma:qIso}
Given chain maps $f\colon C_1\longrightarrow C_2$ and $f'\colon
C_1'\longrightarrow C_2'$ and quasi-isomorphisms $\phi_1\colon
C_1\longrightarrow C_1'$ and $\phi_2\colon C_2\longrightarrow C_2'$ 
so that the composites $f'\circ \phi_1$ is chain homotopic to $\phi_2\circ f$,
there is an induced quasi-isomorphism $\Phi$ from $\MCone(f)$ to $\MCone(f')$
making the diagram
$$
\begin{CD}
0@>>> C_2 @>>> \MCone(f) @>>> C_1 @>>> 0 \\
&& @V{\phi_2}VV @V{\Phi}VV @V{\phi_1}VV \\
0@>>> C_2' @>>> \MCone(f') @>>> C_1' @>>> 0 \\
\end{CD}
$$
commutative.
\end{lemma}

\begin{proof}
Define $\Phi(a_1 \oplus a_2) = (\phi_1(a_1), (-1)^{\deg(a_1)}
H(a_1)+\phi_2(a_2))$, where $H\colon C_1 \longrightarrow C_2'$ is the
homotopy between $f'\circ \phi_1$ and $\phi_2\circ f$. It is easy to
see that $\Phi$ is a chain map, and indeed that the labelled diagram
is commutative. The map $\Phi$ is a quasi-isomorphism by the
five-lemma.
\end{proof}

\subsection{Heegaard Floer complexes}
\label{subsec:CFcx}

A {\em pointed Heegaard diagram for a three-manifold $Y$} is a quadruple
$(\Sigma,\alphas,\betas,z)$, where $\Sigma$ is an oriented surface of
genus $g$, $\alphas=\{\alpha_1,...\alpha_g\}$ and
$\betas=\{\beta_1,...,\beta_g\}$ are complete sets of attaching
circles which specify $Y$, and 
$$z\in
\Sigma-\alpha_1-...-\alpha_g-\beta_1-...-\beta_g$$ is a reference
point. It is explained in~\cite{HolDisk} that this data, together with
some additional analytical choices (including a complex structure over
$\Sigma$) leads to a collection chain complexes
$\CFm(\Sigma,\alphas,\betas,z)$, $\CFinf(\Sigma,\alphas,\betas,z)$,
$\CFp(\Sigma,\alphas,\betas,z)$, and $\CFa(\Sigma,\alphas,\betas,z)$,
which we refer to simply as $\CFc(\Sigma,\alphas,\betas,z)$.
These complexes are constructed from a suitable
variant of Lagrangian Floer homology in the $g$-fold symmetric
product of $\Sigma$.  Specifically, letting
\begin{eqnarray*}
\Ta=\alpha_1\times...\times\alpha_g\subset \Sym^g(\Sigma) &{\text{and}}&
\Tb=\beta_1\times...\times\beta_g\subset \Sym^g(\Sigma),
\end{eqnarray*}
the complex
$\CFinf(\Sigma,\alphas,\betas,z)$ is freely generated over $\Z$ by
generators $[\x,i]\in (\Ta\cap\Tb)\times \Z$, endowed with a
differential $$\partial[\x,i]=\sum_{\y\in\Ta\cap\Tb}
\sum_{\{\phi\in\pi_2(\x,\y)\big| \Mas(\phi)=1\}}
\#\UnparModFlow(\phi)[\y,i-n_z(\phi)],$$
where $\pi_2(\x,\y)$ denotes the space of homotopy classes of Whitney
disks from $\x$ to $\y$, $\Mas(\phi)$ denotes its Maslov index,
$\UnparModFlow(\phi)$ denotes the moduli space of pseudo-holomorphic
representatives of $\phi$ (with respect to some suitably generic
perturbation), divided out by the natural translation action, and
$n_z(\phi)$ denotes the algebraic intersection number of $\phi$ with
the locus $\{z\}\times \Sym^{g-1}(\Sigma)$.  $\CFm(\Sigma,\alphas,\betas,z)$
is the subcomplex generated by $[\x,i]$ with $i<0$, $\CFp(\Sigma,\alphas,\betas,z)$
is its quotient complex (i.e. generated by $i\geq 0$), 
and $\CFa(\Sigma,\alphas,\betas,z)$ is the subcomplex of $\CFp$
generated by pairs with $i=0$. These complexes
are modules over the polynomial algebra $\Z[U]$, where $U\cdot [\x,i]=[\x,i-1]$.
We have induced $\Z[U]$-actions on their homology groups
$\HFinf(Y)$, $\HFm(Y)$, $\HFp(Y)$, and $\HFa(Y)$ respectively.

The homology groups of these complexes are the Heegaard Floer homology
groups of $Y$. As the notation suggests, although the chain complexes
are defined using several choices (including a Heegaard diagram for
$Y$ and a choice of complex structure over $\Sigma$), the homology
groups on only the homeomorphism type of the underlying
three-manifold. In fact, the proof of topological invariance
from~\cite{HolDisk} actually proves more; it shows that
that the chain homotopy type
of the chain complexes $\CFc(\Sigma,\alphas,\betas,z)$ is a
topological invariant of $Y$.  As a shorthand, we let $\CFc(Y)$ denote
the Heegaard Floer complex of $Y$ for some choice of Heegaard diagram
(and auxilliary choices).

The Heegaard Floer complexes come with some additional structure. For
example, there is a splitting of $\CFp(Y)$ into summands indexed by
$\SpinC$ structures over $Y$, $$\CFp(Y)=\bigoplus_{\spinct\in\SpinC(Y)}
\CFp(Y,\spinct).$$ These complexes are typically $\Zmod{2}$-graded,
but when $Y$ is a rational homology three-sphere, or more generally,
when $\spinct$ is a $\SpinC$ structure whose first Chern class is
torsion, then $\CFp(Y,\spinct)$ is naturally a $\Q$-graded complex
(cf.~\cite{HolDiskFour}, see also~\cite{AbsGraded}).

In a slight abuse of notation, we will in fact write $\CFp(Y)$ for a
chain complex which is quasi-isomorphic to a chain complex for $Y$
with respect to some choice of Heegaard diagram. In the cases where
$\CFp(Y)$ has extra structure (e.g. in the case $Y$ is a rational
homology sphere and hence its Heegaard Floer compexes are
$\Q$-graded), we require that our candidate have that additional
structure, and the quasi-isomorphism preserves it.

It is shown in~\cite{HolDiskFour} that Heegaard Floer homology is
natural under cobordisms. Specifically, in the present paper,
three-manifolds $Y$ will always be closed and oriented. A cobordism
$W$ from $Y_1$ to $Y_2$ is a smooth, connected, compact four-manifold
with two boundary components $-Y_1$ and $Y_2$ (with respect to their
boundary orientations). We sometimes write this as $W\colon Y_1
\longrightarrow Y_2$.  Of course, such a cobordism can be ``turned
around'' and viewed as a cobordism $W\colon -Y_2 \longrightarrow
-Y_1$. Given a cobordism $W\colon Y_1 \longrightarrow Y_2$, we say
that two $\SpinC$ structures $\spinct_i\in\SpinC(Y_i)$ for $i=1,2$ are
{\em $\SpinC$ cobordant} if there is a $\SpinC$ structure
$\spinc\in\SpinC(W)$ with $\spinc|Y_{i}=\spinct_i$ for $i=1,2$.

Suppose now that $W\colon Y_1\longrightarrow Y_2$ is a cobordism, equipped
with a $\SpinC$ structure $\spinc$ whose
restrictions $\spinct_i=\spinc|_{Y_i}$ for $i=1,2$ both have  torsion first
Chern class, then there is an induced chain map $$\fp_{W,\spinc}\colon
\CFp(Y_1,\spinct_1)
\longrightarrow \CFp(Y_2,\spinct_2)$$ which is homogeneous
of degree 
\begin{equation}
\label{eq:DimensionShift}
\frac{c_1(\spinc)^2-2\chi(W)-3\sigma(W)}{4}
\end{equation}
(cf. Theorem~\ref{HolDiskFour:thm:AbsGrade} of~\cite{HolDiskFour}).
In fact, the $\Q$-grading on Floer homology is characterized by the
above formula, and the normalization that $\HFp_d(S^3)$ is trivial for
all $d<0$, non-trivial in degree $d=0$.

Sometimes, we will find it convenient to pass to a variant of
Heegaard Floer homology parameterized by an integer $\delta\geq 0$
which interpolates between $\CFa$ and $\CFp$, which we write as
$\CFd$. The generators of $\CFd$, now, are pairs $[\x,i]$ where $0\leq
i \leq \delta$, endowed with the induced differential from $\CFp$. In
other words, $\CFd$ is the subcomplex of $\CFp(Y)$ which is the kernel
of multiplication by $U^{\delta+1}$ (in particular, in the case where
$\delta=0$, this construction gives $\CFa$).

\subsection{Integral surgeries on knots}
\label{subsec:IntegralSurgeries}

Let $K\subset Y$ be a null-homologous knot. Then, there is a canonical
Seifert framing on $K$, giving rise to a curve $\lambda$ in the
boundary of the tubular neighborhood of the knot $K$, $\nbd{K}$, which
meets the meridian $\mu$ in a single point. Given any integer $n$, the
three-manifold $Y_{n}(K)$ denotes the new three-manifold obtained by
$n$-Dehn filling on the complement $Y-\nbd{K}$; i.e. this is a
three-manifold obtained by filling $Y-\nbd{K}$ with a solid torus
whose new meridian is given by $n\cm \mu+\lambda$. Of course, Dehn
filling makes sense for arbitrary rational numbers, but we restrict
attention here to integral surgeries.  For such a surgery, in fact,
there is also a canonical four-manifold $W_{n}(K)$ which is obtained
by attaching a two-handle to $[0,1]\times Y$ with framing $n$ along
$K$. This gives a cobordism from $Y$ to $Y_n(K)$. In our applications,
we find it sometimes convenient to consider the cobordism
$W'_n(K)\colon Y_{n}(K) \longrightarrow Y$ obtained by turning around
the cobordism $-W_n(K)\colon -Y \longrightarrow -Y_n(K)$.

Fix a Seifert surface $F$ for $K$, and let $\CapSurf\subset W'_n(K)$
denote the surface obtained by capping off $F$ in $W'_n(K)$.
Suppose that $\spincu$ is a $\SpinC$ structure over $Y_n(K)$
which admits an extension $\spinc$ over $W'_n(K)$ with the
property that
$$\langle c_1(\spinc),[{\widehat F}]\rangle -n \equiv 2i\pmod{2n}.$$

\begin{lemma}
\label{lemma:IdentifySpinC}
The correspondence $\spincu\mapsto i$ determined by the above
formula induces a surjection $\SpinC(Y_n(K)) \longrightarrow \Zmod{n}$.
In the case where $Y$ is an integer homology three-sphere, the map
is an isomorphism. 
More generally, if we fix a $\SpinC$ structure
$\spinct$ over $Y$, the set of $\SpinC$ structures over $Y_n(K)$ which
are $\SpinC$-cobordant to $\spinct$ over $W'_n(K)$ is identified with $\Zmod{n}$ under this correspondence.
\end{lemma}

\begin{proof}
Straightforward.
\end{proof}

Of course, in the case where $Y$ is an integral homology three-sphere,
the above correspondence between $\SpinC$ structures over $Y_n(K)$ and
$\Zmod{n}$ depends on the Seifert surface only through its induced
orientation on $K$. Moreover, if we fix $i\in\Zmod{n}$ and use the
opposite orientation on $K$, the induced $\SpinC$ structure over
$Y_n(K)$ is conjugated. 

For $i\in \Zmod{n}$, when $Y$ is an integral homology three-sphere we
write $\HFp(Y_n(K),i)$ for the Heegaard Floer homology of $Y_n(K)$
calculated in the $\SpinC$ structure corresponding to $i\in \Zmod{n}$
under the correspondence from the above lemma.

For $n\neq 0$, the lens space $L(n,1)$ can be viewed as $n$-surgery on the unknot.
For $n>0$, let
\begin{equation}
\label{eq:DefOfD}
d(n,i) = -\max_{\{s\in\Z\big|s \equiv i\pmod{n}\}}
\frac{1}{4}\left(1-\frac{(n+2s)^2}{n}\right)
\end{equation}
and let $d(-n,i)=-d(n,i)$.
It can be shown that
$d(n,i)$ is the smallest degree in which $\HFp(L(n,1),i)$ is
non-trivial (cf. Proposition~\ref{AbsGraded:prop:dLens} of~\cite{AbsGraded}).

\subsection{Knot Floer homology and large $n$ surgeries}
\label{subsec:LargeNSurgeries}

We follow here the notation on knot Floer homology from~\cite{Knots}).
A knot $K\subset Y$ has a compatible Heegaard diagram
$(\Sigma,\alphas,\betas,w,z)$, where $(\Sigma,\alphas,\betas)$ is a
Heegaard diagram for $Y$, the knot $K$ is supported in the handlebody
specified by $\betas$, where it is a standard unknotted circle which
is dual to the $\beta_g$-attaching disk, and $w$ and $z$ are a pair of
reference points close to, and lying on either side of $\beta_g$.

This gives rise to a map $$\sRelSpinC \colon \Ta\cap \Tb\longrightarrow
\Z$$ which is half of the first Chern class of the relative $\SpinC$ structure
belonging to $\x\in \Ta\cap\Tb$ evaluated on a Seifert surface for
$K$.  The knot chain complex $C$ described in the introduction is then
generated by $[\x,i,j]\in (\Ta\cap\Tb)\times \Z\times \Z$, satisfying
$\sRelSpinC(\x)+(i-j)=0$, endowed with the differential $$\partial
[\x,i,j] = \sum_{\y\in\Ta\cap\Tb}
\sum_{\{\phi\in\pi_2(\x,\y)\big|\Mas(\phi)=1\}}
\#\UnparModFlow(\phi) \cm [\y,i-n_w(\phi),j-n_z(\phi)].$$
This complex is given the filtration function
$\Filt[\x,i,j]=(i,j)$. The forgetful map $[\x,i,j]\longrightarrow
[\x,i]$ induces an isomorphism between $C$ and $\CFinf(Y)$,
sending $C\{i\geq 0\}$ isomorphically to $\CFp(Y)$.

The following result is proved
in Theorem~\ref{Knots:thm:LargePosSurgeries} of~\cite{Knots}, but we
sketch the proof again for the reader's convenience. (In
fact, a corresponding statement holds for a null-homologous knot in an
arbitrary closed, oriented three-manifold.)

\begin{theorem}
\label{thm:LargeNSurgery}
Let $K\subset Y$ be a null-homologous knot in an integral homology
three-sphere.  There is an integer $N$ with the property that for all
$m\geq N$ and all $t\in\Zmod{m}$, $\CFp(Y_m(K),t)$ is represented by
the chain complex $\Ap_s = C\{\min(i,j-s)\geq 0\}$ where $s\equiv
t\pmod{m}$ and $|s|\leq m/2$, in the sense that there are isomorphisms
(of relatively $\Z$-graded $\Z[U]$-complexes) $$\Psi^+_{m,s}\colon
\CFp(Y_m(K),t)\longrightarrow \Ap_s.$$ Moreover, if ${\mathfrak x}_s$
and ${\mathfrak y}_s$ denote the $\SpinC$ structures over $W_m'(K)$
with
\begin{eqnarray*}
\langle c_1(\spincx_s),[\CapSurf]\rangle +m = 2s
&{\text{resp.}}&
\langle
c_1(\spincy_s),[\CapSurf]\rangle -m = 2s,
\end{eqnarray*}
then $\vertp_s$ and $\horp_s$ correspond to the maps induced by the
cobordism $W'_m(K)$ endowed with the $\spincx_s$ and $\spincy_s$
respectively, in the sense that the following squares commute:
\begin{eqnarray*}
\begin{CD}
\CFp(Y_m(K),t) @>{\fp_{W'_m(K),\spincx_s}}>> \CFp(Y) \\
@V{\Psi^+_{m,s}}VV @VV{=}V \\
\Ap_s @>{\vertp}>> \Bp 
\end{CD}
&{\text{and}}&
\begin{CD}
\CFp(Y_m(K),t) @>{\fp_{W'_m(K),\spincy_s}}>> \CFp(Y) \\
@V{\Psi^+_{m,s}}VV @VV{=}V \\
\Ap_s @>{\horp}>> \Bp.
\end{CD}
\end{eqnarray*}
\end{theorem}

\vskip.2cm
\noindent{\bf{Sketch of Proof.}}
Let $(\Sigma,\alphas,\gammas,\betas,w,z)$ be a Heegaard diagram
for the cobordism $W'_m\colon Y_m(K) \longrightarrow Y$, containing
the pair of basepoints $w$ and $z$ one on each side of the meridian for $Y$.
In particular, here the three-manifold $Y_{\alpha,\gamma}\cong Y_m(K)$,
$Y_{\gamma,\beta}\cong \#^{g-1}(S^2\times S^1)$, $Y_{\alpha,\beta}\cong Y$.

For $s\in\Z$, we define the map $\Psi^+_{m,s}\colon \CFp(Y_m(K),[s])
\longrightarrow C\{\max(i,j-s)\geq 0\}$ by $$\Psi^+_{m,s}[\x,i] =
\sum_{\y\in\Ta\cap\Tb}
\sum_{\{\psi\in\pi_2(\x,\Theta,\y)\big|
n_w(\psi)-n_z(\psi)=s-\sRelSpinC(\x)\}}
\#\ModFlow(\psi)[\y,i-n_w(\psi),i-n_z(\psi)],$$
where here $\Theta\in \HFleq(Y_{\beta,\gamma})$ is a generator.
The constraint on homotopy classes of triangles is equivalent to the constraint
that the first Chern class $\spinc$ induced by the homotopy class
$\psi\in\pi_2(\x,\Theta,\y)$ and the basepoint $w$
satisfies
$$\langle c_1(\spinc),[\CapSurf] \rangle + m = 2s.$$

It is shown in the proof of \cite[Theorem~\ref{Knots:thm:LargePosSurgeries}]{Knots}
for all sufficiently large $m$ and all $t\in \Zmod{m}$, we have
that $\Psi^+_{m,s}$ induces an isomorphism of chain complexes from
$\CFp(Y_m(K),t)$ to $C\{\max(i,j-s)\geq 0\}$ for 
$|s|\leq \frac{m}{2}$ and $s\equiv t\pmod{m}$.
Post-composing $\Psi^+_{m,s}$ with the projection
$C\{\max(i,j-s)\geq 0\}$ to $C\{i\geq 0\}$ (i.e. $\vertp$), we obtain
the map induced by the cobordism $W'_{m}(K)$ equipped with the
$\SpinC$ structure $\spincx_s$; i.e.  the square in the statement of
the theorem on the left commutes. 

For commutativity of the second square, observe that the composite
$\horp\circ \Psi^+_{m,s}$ is identified with a map $\CFp(Y_m(K))$ to
$\CFp(Y)$ induced by counting pseudo-holomorphic triangles in a fixed
$\SpinC$ equivalence class, using the reference point $z$ (rather than $w$).
As we have seen, with respect to the reference point $w$, this $\SpinC$ equivalence
class induces the $\SpinC$ structure $\spincx_s$; thus, with respect to the reference 
point $z$, the induced $\SpinC$ structure is $\spincy_s=\spincx_s -\PD[\CapSurf]$.
\qed
\vskip.2cm

It follows from the above statement that there are corresponding
identifications between $\CFa(Y_m(K),i)$ with $C\{\max(i,j-s)=0\}$, as
well.  

\begin{cor}
  \label{cor:LargeNFloerHomology}
  Let $K\subset Y$ be a null-homologous knot in an integer homology
  three-sphere, and fix an integer $\delta\geq 0$. There are constants
  $C_1$ and $C_2$ (depending on the knot $K$ and the choice of $\delta\geq 0$)
  with the property that for all sufficiently large
  $N$, and any choice of $i$,
  there is a chain complex $\CFd(Y_N,i)$ such that
\begin{eqnarray}
\max\gr\CFd(Y_N,i)-\min\gr\CFd(Y_N,i)&\leq& C_1 
\label{ineq:SpreadYN}\\
\Big|\max\gr\CFd(Y_N,i)-\min\gr\CFd(Y) + 
\max_{s\equiv i\pmod{N}} \frac{1}{4}\Big(1-\frac{|2s+N|^2}{N}\Big)\Big|
&\leq & C_2.
\label{ineq:SpreadYYN}
\end{eqnarray}
\end{cor}

\begin{proof}
  Both statements follow from Theorem~\ref{thm:LargeNSurgery}.  The
  first statement is an immediate consequence of the homogeneous
  identification for all sufficiently large $N$ of $\CFp(Y_N(K),i)$
  with $\Ap_i$ (which in turn is independent of $N$). For the second
  assertion, observe that either $\spincx_i$ or $\spincy_i$ is the
  $\SpinC$ structure over $W'_{N}(K)$ for which $c_1(\spinc)^2$ is
  maximized, amongst all $\SpinC$ structures whose restriction to
  $Y_N(K)$ corresponds to $i$. (In fact, it is $\spincx_i$ if $i\geq
  0$, and $\spincy_i$ if $i\leq 0$.) Morever, the induced map on
  $W'_N(K)$ is realized as $\vertp$ or $\horp$, according to whether
  or not $i\geq 0$. It is now immediate to see that for $i\geq 0$ resp
  $i\leq 0$, there is a constant $C_2$ with the property that all of
  $\CFp(Y)$ is contained within $C_2$ grading levels of the image of
  $\vertp$ resp. $\horp$. The second inequality follows from this
  observation, together with the is a shift in the gradings in
  Equation~\eqref{eq:DimensionShift} specialized to the cobordism
  $W_N'(K)$. 
\end{proof}

\subsection{An example: the unknot} 
\label{subsec:Unknot}
In this subsection, we consider a special case of
Theorem~\ref{thm:IntegerSurgeries}: the case of the unknot $K$ in the
three-sphere.  

For the unknot $K$, the knot Floer complex $\CFKinf(K)$ is generated
over $\Z$ by a sequence of generators $\{x_i\}_{i\in\Z}$ with
$\Filt(x_i)=(i,i)$, and $U(x_i)=x_{i-1}$.

Let $\InjMod{}$ denote the $\Z[U]$ module $\Z[U,U^{-1}]/U\cm\Z[U]$.

\begin{figure}
\mbox{\vbox{\epsfbox{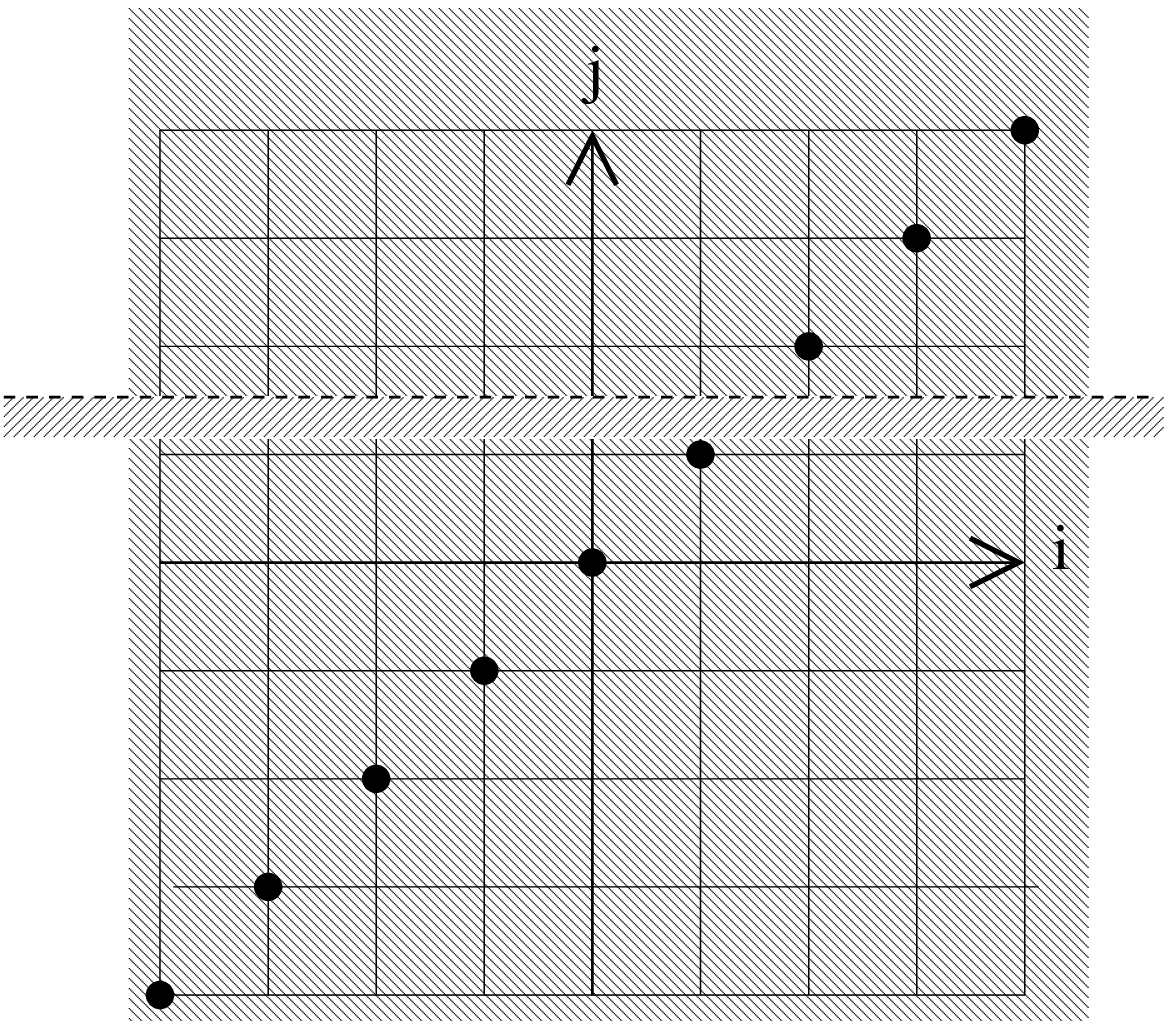}}}
\caption{\label{fig:Unknot}
{\bf The doubly-filtered knot complex for the unknot.}
We have illustrated the doubly-filtered complex for the unknot:
each dot at the location $(i,j)$ 
denotes a $\Z$-generator for the complex $C$ whose filtration level
is $(i,j)$ (and indeed, they
are all located at $(i,i)$). The hatched region represents $\Ap_2$, while projection
to the region above the dotted line represents $\horp_2$. The two generators
in $\Ap_2$ but not above the dotted line represent the two kernel elements
of $\horp_2$.}
\end{figure}

For the unknot, we have that $\Ap_s=\InjMod{}$ for all
$s$, and also $\Bp_s=\InjMod{}$. Moreover, the maps
\begin{eqnarray*}
\vertp_s\colon \Ap_s \longrightarrow \Bp_s
&{\text{and}}&
\horp_s\colon \Ap_s \longrightarrow \Bp_{s+n}
\end{eqnarray*}
can be explicitly identified (up to multiplication by $\pm 1$)
as endomorphisms of $\InjMod{}$:
\begin{eqnarray*}
\vertp_s = \left\{\begin{array}{ll}
1 & {\text{if $s\geq 0$}} \\
U^{-s} & {\text{if $s\leq 0$}}\end{array}\right. \\
\horp_s = \left\{\begin{array}{ll}
U^{s} & {\text{if $s\geq 0$}} \\
1 & {\text{if $s\leq 0$.}} 
\end{array}\right.
\end{eqnarray*}

For fixed $i\in\Zmod{n}$, let $\BigAp_i=\bigoplus_{s\equiv i\pmod{n}} \Ap_s$ and 
$\BigBp_i=\bigoplus_{s\equiv i\pmod{n}} \Bp_s$, and 
$$\Dp_{n,i}\colon \BigAp_i\longrightarrow \BigBp_i$$
denote the map obtained by restricting $\Dp_n$.
Clearly, the mapping cone of $\Dp_n$ splits as a direct sum of the 
mapping cones $\Dp_{n,i}$ over all $i\in\Zmod{n}$.

It is easy to see that the homology of the mapping cone
$$\Dp_{n,i}\colon \BigAp \longrightarrow \BigBp$$ is isomorphic to
$\InjMod{}$. This, of course, is consistent with
Theorem~\ref{thm:IntegerSurgeries}, together with the fact that the
Floer homology of the lens space $L(n,1)$ (which in turn is $n$
surgery on the unknot) has $\HFp(L(n,1),i) \cong \InjMod{}$ for each
$i\in\Zmod{n}$ (cf. Proposition~\ref{HolDisk:prop:Lensspaces}
of~\cite{HolDisk}).

More explicitly, in the case where $n>0$,
it is easy to see that
$\Dp_{n,i}$ is surjective, 
and its kernel is identified with the image of an injection
$$\iota\colon \InjMod{} \longrightarrow \BigAp$$
given by the formula
\begin{equation}
\label{eq:DefIota}
\iota(\xi) = \{U^{\epsilon(\sigma+kn,n)}\cm \xi\}_{k\in\Z},
\end{equation}
where here
$\sigma$ is the representative for $i\pmod{n}$ with $0\leq \sigma <n$, and
\begin{equation}
\epsilon(\sigma+kn,n) = 
\left\{\begin{array}{ll}
k\sigma + \left(\frac{k(k-1)n}{2}\right) & {\text{if $k\geq 0$}} \\
(k+1)\sigma + \left(\frac{k(k+1)n}{2}\right) & {\text{if $k<0$.}}
\end{array}
\right.
\end{equation}
Indeed, $\iota$ induces the isomorphism of $\InjMod{}$ with
$H_*(\MCone(\Dp_{n,i}\colon \BigAp \longrightarrow \BigBp))$.

Dually, if $n>0$, then $\Dp_{-n,i}$ is injective, and its cokernel
is isomorphic to $\InjMod{}$. In fact, for any $\delta\geq 0$,
let $\BigAd$ resp. $\BigBd$ denote the kernel of $U^{\delta+1}$ on 
$\BigAp$ resp $\BigBp$, and let
$$\Dd_{-n,i}\colon \BigAd \longrightarrow \BigBd$$
denote the restriction of $\Dp_{-n,i}$. There is a
map 
$$\pi\colon \BigBd \longrightarrow \Z[U]/U^{\delta+1}$$
defined by
$$\pi(\{\eta_{\sigma+kn}\}_{k\in\Z})=\sum_{k\in\Z} (-1)^k
U^{\epsilon(-(\sigma+kn),n)} \eta_{\sigma+kn}$$
which vanishes on the image
of $\Dd_{-n,i}(\Ad)$, and inducing an isomorphism
$$H_*(\MCone(\Dd_{-n,i}\colon \BigAd\longrightarrow \BigBd)\cong \Z[U]/U^{\delta+1}.$$

\subsection{Algebra}
\label{subsec:Algebra}

Consider the polynomial algebra $\Z[U]$ in a single variable. A 
{\em graded $\Z[U]$ module} is a $\Q$-graded module $M$ over 
the ring $\Z[U]$ with the property that the endomorphism
$U\colon M \longrightarrow M$ is homogeneous with degree $-2$.
For $d\in \Z$, let $M_d$ denote the the subgroup generated by homogeneous
elements of degree $d$, and let $M_{\leq k}$ denote the sum
$$M_{\leq k} =\bigoplus_{\{d\leq k\}} M_d.$$ Note that for all integers
$k$, $M_{\leq k}$ is a $\Z[U]$-submodule of $M$.

We say that a graded $\Z[U]$-module is {\em of $\HFp$ type} if
for all sufficiently large degrees $k$, the endomorphism
$$U\colon M_{d}\longrightarrow M_{d-2}$$ is an isomorphism,
and also for all sufficiently small degrees, $M_d=0$.
Of course, if $Y$ is any three-manifold and $\spinct\in\SpinC(Y)$
is a $\SpinC$ structure whose first Chern class is torsion, then
$\HFp(Y,\spinct)$ is a 
$\Z[U]$-module of $\HFp$-type.

The following lemma is straightforward:

\begin{lemma}
 \label{lemma:GradedModules}
Let $A$ and $B$ be a pair of graded $\Z[U]$-modules of $\HFp$-type.
For any integer $c$ there is a constant $D$ with the property that for
all $k\geq D$, any homogeneous map of degree $c$ defined on $A_{\leq k}$
$$f_{\leq k}\colon A_{\leq k} \longrightarrow B_{\leq c+k}$$
can be uniquely extended to a $\Z[U]$-module map
$$f\colon A \longrightarrow B.$$
In particular, if $A$ and $B$ are of $\HFp$-type then there is a $\delta\geq 0$
with the property that if $A_{\leq \delta}\cong B_{\leq \delta}$, then 
$A\cong B$.
\end{lemma}

\begin{proof}
  Since $A$ and $B$ are of $\HFp$ type, for all sufficiently large
  $\ell$ and all non-negative integers $m$, the maps $$U^{m}\colon
  A_{\ell+2m} \longrightarrow A_{\ell}$$
  and $$U^{m}\colon
  B_{\ell+2m+c} \longrightarrow B_{\ell+c}$$
  are isomorphisms.  Choose
  then $k\geq \ell+2$, and define $$f(a)=\left\{ \begin{array}{ll}
      f_{\leq k}(a) & \text {if $\gr(a)\leq k$} \\ U^{-m} f_{\leq
        k}(U^m(a)) & \text
      {where $m$ is chosen so that $\ell \leq \gr(a)-2m\leq\ell+2$}. \\
  \end{array}\right.  $$ It is easy to see that $f$ is a
  canonically-defined extension of $f_{\leq k}$.  The last claim
  follows immediately.
\end{proof}

\begin{defn}
\label{def:CFpType}
A chain complex $C$ over the ring $\Z[U]$ is said to be {\em of
$\CFp$-type} if it is quasi-isomorphic (over $\Z[U]$) to a chain
complex of the form $C'\otimes_{\Z[U]}
\Z[U,U^{-1}]/\Z[U]$, where $C'$ is a finitely generated, free chain
complex over $\Z[U]$.
\end{defn}

For example, the chain complex $\CFp(Y)$ is of
$\CFp$-type (the chain complex $\CFm(Y)$ plays the role here of $C'$).
It is also straightforward to see that if $C$ is a graded
$\Z[U]$-complex which is of $\CFp$-type, then its homology is of
$\HFp$-type.  The mapping cones $\Xp(n)$ from the introduction can
also be seen to be of $\CFp$-type (cf. Lemma~\ref{lemma:TruncateX}
below).

Let $C$ be a chain complex over $\Z[U]$, and fix an integer
$\delta\geq 0$.  Let $C^{\delta}$ denote the subcomplex of elements in
the kernel of multiplication by $U^{\delta+1}$.

\begin{lemma}
  \label{lemma:HFdLemma}
  Let $A$ and $B$ be two chain complexes of graded $\Z[U]$-modules
  which are of $\CFp$-type. For any $c$, there is a constant
  $D$ with the property that for all integers
  $\delta \geq D$, any homogeneous map of
  degree $c$ defined on
  $$ 
  f^{\delta}\colon H_*(A^{\delta}) \longrightarrow H_*(B^{\delta})
  $$
  can be uniquely extended to a $\Z[U]$-module map
  $$
  f\colon H_*(A) \longrightarrow H_*(B).
  $$
  In particular, if $A$ and $B$ are graded $\Z[U]$-complexes
  which are of $\CFp$-type, then there is an integer $\delta\geq 0$ with the
  property that if $H_*(A^{\delta}) \cong H_*(B^{\delta})$, then
  $H_*(A)\cong H_*(B)$.
\end{lemma}

\begin{proof}
  Consider the short exact sequence
  $$ 
  \begin{CD}
    0@>>> A^\delta @>>> A@>{U^{\delta+1}}>> A@>>> 0.
  \end{CD} 
  $$
  note that the connecting homomorphism is a homogeneous map of
  degree $2\delta+1$. In particular, for any $k$, we can choose
  $\delta$ so that $H_*(A_{\leq k-2\delta-1})=0$, the induced long exact
  sequence gives an isomorphism $H_*(A^{\delta}_{\leq k})\cong
  H_*(A_{\leq k})$. The lemma then follows at once from
  Lemma~\ref{lemma:GradedModules}.
\end{proof}

\section{An exact sequence for surgeries}
\label{sec:ExactSeq}

Theorem~\ref{thm:IntegerSurgeries} hinges on the following exact
sequence relating different surgeries on a knot in a three-manifold
$Y$.  Loosely speaking, the exact sequence is induced by the homology
between the curve  $(m+n,1)$ in the torus with the sum of curves
$(n,1)$ and $(m,0)$ (where here $m$, $n$, are integers and
the homology classes of curves in the torus are written
as $(a,b)\in\Z\oplus\Z\cong H_1(T^2;\Z)$).

\begin{theorem}
\label{thm:ExactSeq}
Let $Y$ be a closed, oriented three-manifold, equipped with a
null-homologous knot $K$.  Fix an integer $n$ and a positive integer
$m$.  Then, there is a long exact sequence
$$\begin{CD}
...@>>>\HFp(Y_n(K)) @>>>\HFp(Y_{m+n}(K)) @>>> \bigoplus^m \HFp(Y) @>>>...;
\end{CD}$$
and also a corresponding exact sequence using $\HFa$ in place of $\HFp$.
Indeed, there are $\Z[U]$-equivariant chain maps 
\begin{eqnarray*}
\fp_1\colon \CFp(Y_n(K)) &\longrightarrow& \CFp(Y_{m+n}(K)) \\
\fp_2\colon \CFp(Y_{m+n}(K)) &\longrightarrow& \bigoplus^m \CFp(Y) \\
\fp_3\colon \bigoplus^m \CFp(Y) &\longrightarrow& \CFp(Y_n(K))
\end{eqnarray*}
inducing the maps in the long exact sequence,
and $\Z[U]$-equivariant quasi-isomorphisms
\begin{eqnarray*}
\phi^+\colon \CFp(Y_n(K))&\longrightarrow & \MCone(\fp_2) \\
\psi^+\colon \MCone(\fp_2) &\longrightarrow & \CFp(Y_n(K)).
\end{eqnarray*}
(See Equations~\eqref{eq:DefF1}, \eqref{eq:DefF2}, \eqref{eq:DefF3},
\eqref{eq:DefPhi}, and \eqref{eq:DefPsi}
for the definitions
of $\fp_1$, $\fp_2$, $\fp_3$, $\phi^+$, and $\psi^+$ respectively.)
\end{theorem}

The proof follows very closely along the lines of various other
previously-established exact sequences for surgeries. In 
particular, we assume familiarity with
Section~\ref{HolDiskTwo:sec:Surgeries} of~\cite{HolDiskTwo},
and continue in the notation set up there.

\vskip.2cm
\noindent{\bf{Proof of Theorem~\ref{thm:ExactSeq}.}}
Consider a pointed Heegaard diagram for $Y$
$(\Sigma,\alphas,\betas,z)$, with the property that $K$ is contained
entirely inside the handlebody $U_\beta$, so that it is disjoint from
the attaching disks bounding $\beta_1,...,\beta_{g-1}$, and with the
additional property that $\beta_{g}$ is a meridian for the knot.  Let
$\gamma_{g}$ be a simple, closed curve in $\Sigma$ disjoint from the
$\beta_1,...,\beta_{g-1}$ which specifies the $n$-framing of $K$.  (In
particular, if $\lambda_g$ is the canonical $0$-framing, then
$\gamma_g$ is a smooth curve which is homologous to $n\cm \beta_g +
\lambda_g$.) We complete this to a $g$-tuple of attaching circles
$\gammas$ by taking curves $\gamma_1,...,\gamma_{g-1}$ which are small
Hamiltonian translates of $\beta_1,...,\beta_{g-1}$
respectively. Similarly, define $\deltas$, only this time $\delta_{g}$
corresponds to the framing $m+n$.  Thus, $(\Sigma,\alphas,\gammas)$
and $(\Sigma,\alphas,\deltas)$ are Heegaard diagrams for $Y_n(K)$ and
$Y_{m+n}(K)$ respectively.

We place a basepoint $p$ on $\beta_g$, and
consider twisted homology with coefficients in $\Zmod{m}$; i.e. write
$\Z[\Zmod{m}]=\Z[T]/(T^m-1)$, and consider the chain complex
$\CFp(Y)\otimes_\Z \Z[\Zmod{m}]$
endowed with the differential $${\underline \partial}^+ [\x,i] =
\sum_{\y\in\Ta\cap\Tb}\sum_{\{\phi\in\pi_2(\x,\y)\big|\Mas(\phi)=1
\}}\#\left(\frac{\ModFlow(\phi)}{\R}\right) \cm T^{m_{p}(\phi)}\cm
[\y,i-n_{z}(\phi)]$$ 
where as usual here $\x\in\Ta\cap\Tb$, $i\geq 0$, $\pi_2(\x,\y)$
denotes the space of homotopy classes of Whitney disks 
connecting $\x$ and $\y$, $\Mas(\phi)$ denotes the Maslov index  of $\phi$,
and terms in the above equation for which $i-n_z(\phi)<0$ are to be dropped.
Moreover, $m_p(\phi)$ denotes the multiplicity of
the basepoint $p$ in the boundary of $\phi$; i.e. $p$ determines a
codimension one submanifold $\beta_1\times...\times \beta_{g-1}
\times\{p\}\subset \Tb$, and we look at the intersection number
with the restriction of the boundary of $\phi$ with this subset. We
denote the complex by $\CFp(Y;\Z[\Zmod{m}])$ (In the terminology
of~\cite{HolDisk}, this is the chain complex for $Y$ with twisted
coefficients in $\Z[\Zmod{m}]$, where it is denoted
$\uCFp(Y;\Z[\Zmod{m}])$, however, we drop the underline here in the
interest of notational simplicity.) 
There is an isomorphism of chain complexes of modules over
$\Z[\Zmod{m}]$, 
\begin{equation}
\label{eq:DefChangeP}
\theta\colon \CFp(Y;\Z[\Zmod{m}])\stackrel{\cong}{\longrightarrow}
\CFp(Y)\otimes_\Z \Z[\Zmod{m}],
\end{equation}
where here the right-hand-side is
endowed with the differential which is the original differential on
$\CFp(Y)$ tensored with the identity map on $\Z[\Zmod{m}]$.  This map
is induced by fixing an intersection point $\x_0\in \Ta\cap\Tb$, and
defining $$\theta [\x,i]= T^{m_p(\phi)} [\x,i],$$ where here
$\phi\in\pi_2(\x,\x_0)$. Note that $m_p(\phi)$ depends on $\phi$ only
through the choice of $\x_0$ and $\x$ (and indeed the map $\theta$
depends on the choice of $\x_0$ and the placement of $p$ through an overall multiple
of a power of $T$).
There is
a corresponding identification $$\HFp(Y;\Z[\Zmod{m}])\cong
\HFp(Y)\otimes_\Z\Z[\Zmod{m}]\cong \bigoplus^m\HFp(Y).$$

First, we define the map $\fp_1$.  The map $\fp_1$ is defined by
counting pseudo-holomorphic triangles between $\Ta$, $\Tc$, and
$\Td$. More precisely, note that the Heegaard triple
$(\Sigma,\alphas,\gammas,\deltas,z)$ determines a four-manifold
$X_{\alpha,\gamma,\delta}$ with three boundary components
$$Y_{\alpha,\gamma}\cong Y_n(K), \hskip1cm
        Y_{\alpha,\delta}\cong Y_{n+m}(K), 
\hskip1cm \text{and} \hskip1cm Y_{\gamma,\delta}\cong \#^{g-1}(S^2\times S^1)\# L(m,1).$$
We will fix a $\SpinC$ structure over
$Y_{\gamma,\delta}$ by the following convention:

\begin{defn}
  \label{eq:CanonicalSpinC}
  The lens space $L(m,1)$ bounds a tubular neighborhood of a sphere
  $S$ whose self-intersection number is $m$. The {\em canonical $\SpinC$
    structure} $\spinccan\in\SpinC(L(m,1))$ is the one which bounds a
  $\SpinC$ structure $\spinc$ over the tubular neighborhood which
  satisfies $\langle c_1(\spinc),[S]\rangle =m$. For the connected sum
  $L(m,1)\#(\#^{g-1}(S^2\times S^1))$, the canonical $\SpinC$
  structure is the one whose first Chern class is torsion and whose
  restriction to $L(m,1)$ is the canonical $\SpinC$ structure.
\end{defn}

Let $\Theta_{\gamma\delta}$ denote the Floer homology class
corresponding to the generator (over $\Wedge^*
H_1(Y_{\gamma,\delta})\otimes \Z[U]$) of
$$\HFleq(Y_{\gamma,\delta},\spinccan) \cong \Wedge^*
H^1(Y_{\gamma,\delta})\otimes \Z[U]$$ in its canonical $\SpinC$
structure $\spinccan$.  For simplicity, we can arrange for
the homology class
$\Theta_{\gamma\delta}$ to be represented by a single intersection
point in $\Tc\cap\Td$, which we also denote by
$\Theta_{\gamma\delta}$.

We then define
\begin{equation}
\label{eq:DefF1}
f^+_1([\x,i]) = \sum_{\y\in\Ta\cap\Td}\sum_{\{\psi\in\pi_2(\x,\Theta_{\gamma\delta},\y)\big|\Mas(\psi)=0\}}
\#\ModFlow(\psi)\cm [\y,i-n_z(\psi)].
\end{equation}
Similarly, we define $\fp_2\colon \CFp(Y_{m+n}(K))\longrightarrow \CFp(Y;\Z[\Zmod{m}])$ by
\begin{equation}
\label{eq:DefF2}
f^+_2([\y,i]) = \sum_{\w\in\Ta\cap\Tb}\sum_{\{\psi\in\pi_2(\y,\Theta_{\delta\beta},\w)
\big|\Mas(\psi)=0\}} \#\ModFlow(\psi)\cm [\w,i-n_z(\psi)]\cm T^{m_p(\psi)}.
\end{equation}
The map $f^+_2$ can be thought of more invariantly as a weighted sum
of maps induced by cobordisms, as follows.  For our Heegaard triple
$(\Sigma,\alphas,\deltas,\betas,z)$, the doubly-periodic domains all
have multiplicity zero at $p$, while the generator of the space of
triply-periodic domains modulo doubly-periodic ones contains
$\beta_{g}$ with multiplicity $m+n$.  Thus, fixing 
$\y_0$, $\Theta_{\delta\beta}$, and $\w_0$ in $\Ta\cap \Td$,
$\Td\cap \Tb$ and $\Ta\cap \Tb$ respectively, the function
$$m_p\colon \pi_2(\y_0,\Theta_{\delta\beta},\w_0) \longrightarrow
\Zmod{m}$$ descends to a function $$\imult \colon
\SpinC(W'_{n+m}(K))
\longrightarrow \Z$$ with the property that 
$$\imult(\spinc-\PD[{\widehat F}]) = \imult(\spinc)+m+n.$$ 
Note that
$\imult$ depends on the choice of $p$ and the initial fixed intersection points
up to an overall additive
constant; there is also an orientation issue here which is a matter of
convention. 
Composing with the isomorphism from Equation~\eqref{eq:DefChangeP},
it follows at once that
\begin{equation}
\label{eq:InvariantF2}
\theta\circ f^+_2 = \sum_{\spinc\in\SpinC(W'_{n+m}(K))}
T^{\imult(\spinc)}\cm \fp_{W'_{n+m}(K),\spinc}.
\end{equation}

We must check that $f^+_2\circ f^+_1\simeq 0$. This is proved as usual
using associativity for holomorphic triangles. The verification
involves a calculation in the Heegaard triple
$(\Sigma,\gammas,\deltas,\betas,z)$, showing that
the pairing of $\Theta_{\gamma\delta}$ and $\Theta_{\delta\beta}$ in
$\HFleq(Y_{\gamma,\beta};\Z[\Zmod{m}])$ vanishes; i.e.  $$\sum_{\x\in
  \Tc\cap \Tb} \sum_{\{\psi\in\pi_2(\Theta_{\gamma\delta},
  \Theta_{\delta\beta},\x)\big| \Mas(\psi)=0\}} \#\ModFlow(\psi)\cm
[\x,-n_{z}(\psi)]\cm T^{m_p(\psi)}$$
is null-homologous.

The claim amounts to showing that triangles come in cancelling pairs
with the same value of $n_z(\psi)$ and $T$ power. The fact that any
two triangles with the same value of $n_z(\psi)$ connecting
$\Theta_{\gamma\delta}$, $\Theta_{\delta\beta}$, and
$\Theta_{\gamma\beta}$ have the same power of $T$ follows from the
fact that a generator for the space of the triply-periodic domains for
the Heegaard triple $(\Sigma,\gammas,\deltas,\betas,z)$ has
$m_p(\psi)=m$, and of course any two triangles with the same value of
$n_z(\psi)$ differ by a triply-periodic domain.  Next, we claim that
triangles with a fixed value of $n_z(\psi)$ come in pairs. This is
straightforward to verify -- in fact, a linear map transforms the
problem into the same lattice point count
considered in the integer surgeries long exact sequence,
Theorem~\ref{HolDiskTwo:thm:ExactP} of~\cite{HolDiskTwo}.  The usual
proof of the associativity law for maps induced by triangles now gives
a null-homotopy $\Hp{1}\colon \CFp(Y_{n}(K))\longrightarrow
\CFp(Y;\Z[\Zmod{m}])$ of the composite $\fp_2\circ\fp_1$, defined by
counting pseudo-holomorphic quadrilaterals.

One could now establish existence of a long exact sequence in homology
following the outline of~\cite{HolDiskTwo} or~\cite{BrDCov}. Following
the approach of~\cite{HolDiskTwo}, observe that there are curves in
the isotopy class of $\delta_g$ which approximate arbitrarily well the
juxtaposition of curves $m\cm \beta_g+\gamma_g$.  The map from
$\CFp(Y_n(K))$ into $\CFp(Y_{m+n}(K))$ is injective, since the counts
of small triangles (with vertex at the canonical $\SpinC$ structure)
induce the nearest point map from $\CFp(Y_n(K))$ into
$\CFp(Y_{m+n}(K))$, while the counts of small triangles with vertex at
$\beta_g\cap \gamma_g$ and reference point $p\in \beta_g$ induces the
nearest point map tensored with a surjection onto $\Zmod{m}$.

However, for the stronger statement about the quasi-isomorphisms, we
need to use the argument from Section~\ref{BrDCov:sec:LinkSurgeries}
of~\cite{BrDCov}, and in particular, we need to place the third map
on equal footing with the first two. 

To this end, let $c\in\Zmod{m}$ be the element gotten as follows.
Fix any
$\psi\in\pi_2(\Theta_{\gamma\beta},\Theta_{\beta\delta},\Theta_{\gamma\delta})$,
where the $\Theta_{i,j}$ have been chosen before (recall that
$\Theta_{\gamma\delta}$ represents the generator of $\HFleq(L(m,1)\#
(\#^{g-1}(S^2\times S^1)),\spinccan)$ in the canonical $\SpinC$
structure). Observe that the congruence class of
$m_p(\psi)$ modulo $m$ is
independent of the choice of $\psi$ (as any two such $\psi$ differ by
a triply-periodic domain, whose multiplicity at $p$ is some multiple
of $m$). Denote this quantity by $c\in\Zmod{m}$. (It, of course, depends
on the particular Heegaard triple.)

Define
$$\fp_3\colon \CFp(Y;\Z[\Zmod{m}])\longrightarrow \CFp(Y_{n}(K))$$
by the formula
\begin{equation}
\label{eq:DefF3}
\fp_3(T^s\cm [\x,i])= 
\sum_{\y\in\Ta\cap\Tc} \sum_{\{\psi\in\pi_2(\x,\Theta_{\beta\gamma},\y)\big|
\stackrel{\Mas(\psi)=0,}{s+m_p(\psi})\equiv c\pmod{m}
\}}
\#\ModFlow(\psi)\cm [\y,i-n_z(\psi)].
\end{equation}
More invariantly, we can think of $\fp_3$ as a sum of maps induced by 
the two-handle cobordism $W_n(K)\colon Y \longrightarrow Y_n(K)$.
Specifically, multiplicity at $p$ induces a map
$$\imult'\colon \SpinC(W_n(K))\longrightarrow \Z$$
with the property that
$\imult'(\spinc-\PD[{\widehat F}])=\imult'(\spinc)+n$. 
Now, there is some $c'\in\Zmod{m}$ with the property that
\begin{equation}
\label{eq:InvariantF3}
\fp_3\circ \theta(T^s\cm \xi) = \sum_{\{\spinc\in\SpinC(W_n(K))\big| 
\imult'(\spinc)+c'\equiv s\pmod{n}\}} \fp_{W_n(K),\spinc}.
\end{equation}

The verifications that $\fp_3\circ \fp_2\simeq 0$ and $\fp_1\circ
\fp_3\simeq 0$ proceed as before.  In fact, counts of quadrilaterals
define null-homotopies $\Hp{i}$ of $\fp_{i+1}\circ\fp_{i}$. 

For the quasi-isomorphism statement, we count pseudo-holomorphic
quadrilaterals for Heegaard quadruples.  Let $\betas'$, $\gammas'$,
and $\deltas'$ denote suitable exact Hamiltonian translates of the
corresponding curves $\betas$, $\gammas$, and $\deltas$ respectively.
It is easy to see that for the resulting Heegaard tuples
$(\Sigma,\gammas,\deltas,\betas,\gammas',z)$, and
$(\Sigma,\deltas,\betas,\gammas,\deltas',z)$, there is exactly one
Whitney quadrilateral $\varphi$ connecting the $\Theta_{i,j}$ with
$\cald(\varphi)\geq 0$, $n_z(\varphi)=0$, and $\Mas(\phi)=-1$.
(Here, $\cald(\varphi)$ denotes the set of local multiplicities
of $\varphi$ on $\Sigma$, cf.~\cite{HolDisk}; these local multiplicities
must be non-negative for the homotopy class to admit pseudo-holomorphic
representatives). Moreover,
this homotopy class admits a unique pseudo-holomorphic representative, compare
Section~\ref{BrDCov:sec:LinkSurgeries} of~\cite{BrDCov}). 
This is a
routine adaptation of the arguments given in that section.  For the
Heegaard tuple $(\Sigma,\gammas,\deltas,\betas,\gammas',z)$, it
follows that under the map induced by counting holomorphic quadrilaterals, the
map
$$\HFa(\Tc,\Td)\otimes \HFa(\Td,\Tb)\otimes \HFa(\Tb,{\mathbb T}_{\gamma'})
\longrightarrow \HFa(\Tc,{\mathbb T}_{\gamma'})$$
carries the top-dimensional generator of the left-hand-side 
to the top-dimensional generator of the right-hand-side
(cf. \cite[Equation~\eqref{BrDCov:eq:CanonicalSquare}]{BrDCov}).)
A corresponding statement for $\HFleq$ 
(taken with coefficients in the ring of formal power series in $U$,
$\Z[[U]]$) is now a formal consequence,
showing that the top-dimensional generator of the tensor product of 
$\HFleq$ is mapped to a multiple of the
generator of $\HFleq(\Tc,{\mathbb T}_{\gamma'})$
by a unit in $\Z[[U]]$. Corresponding remarks apply to the quadruple
$(\Sigma,\deltas,\betas,\gammas,\deltas',z)$.

The Heegaard tuple $(\Sigma,\betas,\gammas,\deltas,\betas')$ works
slightly differently, as now we are required to take twisted
coefficients when using the pair $\betas$ and $\betas'$,
$\HFa(\Tb,{\mathbb T}_{\beta'},\Z[\Zmod{m}])$. In this case, the
generator of the top-dimensional homology group is not represented by
the intersection point $\Theta_{\beta,\beta'}$, but rather by the
element $(\sum_{i=0}^{m-1} T^i)\cdot \Theta_{\beta,\beta'}$.
Correspondingly, there are now $m$ different homotopy quadrilaterals
$\varphi_i$ with $i=0,...,m-1$ which have $\cald(\varphi_i)\geq 0$,
$n_z(\varphi)=0$, and $\Mas(\varphi)=-1$. We can order these
so that $m_p(\varphi_i)=i$. Each of these has a unique
holomorphic representative.
Again, this shows that the map induced by counting holomorphic
quadrilaterals (and recording their multiplicity at $p$ in the exponent
of a formal variable $T$)
$$\HFa(\Tb,\Tc)\otimes \HFa(\Tc,\Td)\otimes \HFa(\Td,{\mathbb T}_{\beta})
\longrightarrow \HFa(\Tb,{\mathbb T}_{\beta'};\Z[\Zmod{m}])$$
carries the top-dimensional generator of the left-hand-side 
to $(\sum_{i=0}^{m-1}T^i)\cm \Theta_{\beta,\beta'}$, 
the top-dimensional generator of the right-hand-side.
Corresponding remarks for the maps induced on $\HFleq$ follow at once.

With this said, the argument from~\cite{BrDCov} prove that the chain maps
\begin{eqnarray*}
\phi^+\colon \CFp(Y_n(K))\longrightarrow \MCone(\fp_2)
&{\text{and}}&
\psi^+\colon \MCone(\fp_2) \longrightarrow \CFp(Y_{n}(K))
\end{eqnarray*}
defined by 
\begin{equation}
\label{eq:DefPhi}
\phi^+(\xi)=(\fp_1(\xi),\Hp{1}(\xi)).
\end{equation}
and
\begin{equation}
\label{eq:DefPsi}
\psi^+(x,y)=\Hp{2}(x)+\fp_3(y)
\end{equation}
respectively are quasi-isomorphisms.  (Observe, incidentally, that
in~\cite{BrDCov}, we gave a proof using $\CFa$.  To pass from there to
$\CFp$, observe that the chain maps and chain homotopies used are all
$\Z[U]$-equivariant maps; and hence the chain maps above are in fact
$\Z[U]$-equivariant quasi-isomorphisms.)
\qed
 
\section{Proof of Theorem~\ref{thm:IntegerSurgeries}.}
\label{sec:Proof}

Before proving Theorem~\ref{thm:IntegerSurgeries}, we give a more precise
statement.

Let $Y$ be an integral homology three-sphere,
equipped with a knot $K\subset Y$. Fix also an integer $n\neq 0$.

In the introduction, we considered chain complexes
$\BigAp=\bigoplus_{s\in\Z}\Ap_s$ and $\BigBp=\bigoplus_{s\in\Z}\Bp_s$
and a chain map
$\Dp_{n}\colon \BigAp\longrightarrow \BigBp$
defined
by
$\Dp_n(\{a_s\}_{s\in\Z})= \{b_s\}_{s\in\Z},$
where here
$b_s = \horp_{s-n}(a_{s-n})+\vertp_s(a_s)$.
We refine this further as follows.

For each $i\in\Zmod{n}$, let
$\BigAp_i=\bigoplus_{\{s\in\Z\big|s\equiv i\pmod{n}\}} \Ap_s$,
and 
$\BigBp_i=\bigoplus_{\{s\in\Z\big|s\equiv i\pmod{n}\}} \Bp_s$.
There is a splitting
\begin{eqnarray*}
\BigAp=\bigoplus_{i\in\Zmod{n}}\BigAp_i &{\text{and}}&
\BigBp=\bigoplus_{i\in\Zmod{n}}\BigBp_i;
\end{eqnarray*}
and the restriction of $\Dp_n$ to $\BigAp_i$ can be thought of
as a map
\begin{equation}
\label{eq:DefDpN}
\Dp_{n,i}\colon \BigAp_i\longrightarrow \BigBp_i.
\end{equation}

Moreover, the mapping cone $\Xp(n)$ of $\Dp_n$
splits into a direct sum of mapping cones
$\Xp(n)=\bigoplus_{i\in\Zmod{n}} \Xp_{i}(n)$, where here
$\Xp_{i}(n)$ is the mapping cone of $\Dp_{n,i}$.

Clearly, $\Xp(n)$ can be given a relative $\Z$-grading, which
is compatible with the relative $\Z$ grading on $\Ap_s$ and $\Bp_s$,
and with the property that $\vertp$ and $\horp$ (thought of as
endomorphisms of $\Xp(n)$) are homogeneous maps of degree $-1$.

More explicitly, suppose that $n>0$, and consider the absolute grading
on $\CFp(Y)$. Although each $\Bp_s\cong \CFp(Y)$, we shift the
grading; i.e. writing $s=\sigma+\ell\cm n$ where $\sigma$ is the
representative for $i\pmod{n}$ with $0\leq \sigma<n$.  We lift the
natural relative $\Z$-grading on $\Bp_s$ so that the homogeneous
elements of $\CFp(Y)$ of degree $d$ correspond to the homogeneous
elements of $\Bp_{\sigma+\ell\cm n}$ of degree
$d+2\ell\sigma+n\ell(\ell-1)-1$.  We claim we can then consistently
shift the gradings on $\Ap_s$ so that both maps $\vertp_s\colon \Ap_s
\longrightarrow \Bp_s$ and $\horp_s\colon \Ap_s\longrightarrow
\Bp_{s+n}$ are homogeneous maps of degree $-1$.  This gives a grading
on $\Xp_i(n)$.

When considering $\Xp(-n)$ (where here $n>0$), we still have that
$\Bp_s
\cong\CFp(Y)$.  Writing $s=-(\sigma+\ell\cm n)$, for $0\leq
\sigma\leq n$, we identify homogeneous elements of degree $d$ in
$\CFp(Y)$ with the homogeneous elements of $\Bp_{-\sigma-\ell\cm n}$
of degree $d-2\ell \sigma-n\ell(\ell-1)$.
We can then consistently shift the
gradings on $\Ap_s$ so that both maps $\vertp_s\colon \Ap_s
\longrightarrow \Bp_s$ and $\horp_{s}\colon \Ap_s\longrightarrow
\Bp_{s-n}$ are homogeneous maps of degree $-1$.
This gives an integral grading on $\Xp_i(-n)$.

\begin{theorem}
\label{thm:PreciseIntegerSurgeries}
Let $Y$ be an integral homology three-sphere.  The homology of
the mapping cone
$\Xp(n)$ of $$\Dp_n\colon {\mathbb A}^+ \longrightarrow {\mathbb
  B}^+$$
is isomorphic to  $\HFp(Y_n(K))$.
Moreover, the summand $\HFp(Y_n(K),i)$ (under the identification
$\SpinC(Y_n(K))\cong \Zmod{n}$ given in
Subsection~\ref{subsec:IntegralSurgeries}) is identified with the
homology of the
summand $\Xp_i(n)$. Also, this map identifies
the relative $\Z$-grading on $\HFp(Y_n(K),i)$ with the one described
above on $H_*(\Xp_i(n))$.
In fact, the isomorphism from $H_*(\Xp_i(n))$
to  $\HFp(Y_n(K),i)$ is a homogeneous map of degree
$d(n,i)$, as defined in  Equation~\eqref{eq:DefOfD}.
\end{theorem}

The following result is a quick consequence of the proof of
Theorem~\ref{thm:PreciseIntegerSurgeries}. It shows that knot Floer
homology also determines the maps induced by the natural cobordisms
from $Y$ to $Y_n(K)$.

\begin{theorem}
  \label{thm:Maps} The following square commutes: $$ \begin{CD}
  H_*(\BigBp) @>{\iota}>> H_*(\Xp(n)) \\ @V{=}VV @VV{\cong}V \\
  \bigoplus_{i\in\Z}\HFp(Y)_i @>{\Fp{W_n(K)}}>> \HFp(Y_n(K)) \end{CD}
  $$ where the second vertical isomorphism is provided by
  Theorem~\ref{thm:PreciseIntegerSurgeries}, the top horizontal map is
  the natural map induced from the mapping cone construction, while
  the bottom one is the map induced by cobordisms; i.e. the $i^{th}$
  factor $\HFp(Y)_i$ is mapped to $\HFp(Y_n(K))$ by the natural map
  induced by the two-handle cobordism $W_n(K)$ endowed with the
  $\SpinC$ structure with $$\langle c_1(\spinc),[\CapSurf]\rangle
  +n=2i.$$
\end{theorem}

Note that Theorem~\ref{thm:IntegerSurgeries} is a special case of
Theorems~\ref{thm:PreciseIntegerSurgeries} and~\ref{thm:Maps}.

Most of this section is devoted to a proof of
Theorem~\ref{thm:PreciseIntegerSurgeries}. In
Subsection~\ref{subsec:Truncation}, we set
up some notation and algebraic preliminaries for the proof.  The
theorem follows by applying Theorem~\ref{thm:ExactSeq}, choosing the
surgery coefficient $m$ to be sufficiently large. Thus, in in
Subsection~\ref{subsec:GradingsOnW}, we study $\SpinC$ structures over
one of the cobordisms from Theorem~\ref{thm:ExactSeq}, when
the surgery coefficient is large.  This provides enough information to
prove the special case of Theorem~\ref{thm:IntegerSurgeries} for the
case of $\HFa$, cf.  Subsection~\ref{subsec:HFa}.  We assume for
simplicity that $n>0$ throughout most of the exposition, returning to
the case where $n<0$ in a later subsection.  The case for $\HFd$
(arbitrary $\delta\geq 0$) follows as well, with minor notational
modifications, cf.  Subsection~\ref{subsec:CFd}. In
Subsection~\ref{subsec:MoreGradings}, we return to the question of
$\SpinC$ structures over cobordisms, when the surgery coefficient is
sufficiently large. With this done, we can complete the the proof of
Theorem~\ref{thm:PreciseIntegerSurgeries}, in the case where $n>0$. In
Subsection~\ref{subsec:NegativeN}, we describe the modifications
needed to make the proof go over to the case where $n<0$.  In
Subsection~\ref{subsec:ZeroSurgery}, we make a few comments regarding
the case where $n=0$. In Subsection~\ref{subsec:Maps}, we give the
proof of Theorem~\ref{thm:Maps}. In Subsection~\ref{subsec:Generalize},
we discuss some generalizations of Theorem~\ref{thm:PreciseIntegerSurgeries}.

\subsection{Truncation}
\label{subsec:Truncation}

Although we are primarily interested in $\HFp$, we will find it
convenient to use $\HFd$ at various times, and also easier to explain
things in terms of $\HFa$. To this end, we will use some
straightforward notational changes.  Specifically, for any integer
$\delta\geq 0$, let $\Ad_i$, $\Bd_i$, $\BigAd_i$, and $\BigBd_i$
denote the subsets of $\Ap_i$, $\Bp_i$, $\BigAp_i$, and
$\BigBp_i$ respectively which lie in the kernel of multiplication by
$U^{\delta+1}$.  Thus, the mapping cone $\Xd_i(n)$ of the restriction
$\Dd_{n,i}$ of $\Dp_{n,i}$, viewed as a map
$$\Dd_{n,i}\colon \BigAd_i\longrightarrow \BigBd_i,$$
is the kernel of
the action by $U^{\delta+1}$ on $\Xp_i(n)$.  Similarly, we let
$\Aa_i$,$\Ba_i$, $\BigAa_i$, $\BigBa_i$ denote the corresponding groups
$\Ad_i$, $\Bd_i$, $\BigAd_i$, and $\BigBd_i$ with $\delta=0$.

It will be convenient for us also to truncate the infinite
constructions from the introduction.  Specifically, for fixed integers
$\Width$, let
\begin{eqnarray*}
\BigAp_i(\Width)
&=&\bigoplus_{\{s\in\Z\big| s\equiv i\pmod{n}, -\Width\leq s\leq \Width\}} \Ap_s
\subset \BigAp_i \\
\BigBp_i(\Width)
&=&\bigoplus_{\{s\in\Z\big| s\equiv i\pmod{n}, -\Width+n\leq s\leq \Width\}} \Bp_s
\subset \BigBp_i.
\end{eqnarray*}
Similarly, we let
$\Dp_{n,i;\Width}\colon \BigAp(\Width) \longrightarrow \BigBp(\Width)$ 
denote the restriction of $\Dp_{n,i}$, and let
and
$\Xp_i(n;\Width)\subset \Xp_i(n)$ denote the corresponding subset of the
mapping cone. 

Note that when $n>0$, there is a quotient map $Q^+_{\Width}\colon
\Xp_i(n)\longrightarrow \Xp_i(n;\Width)$; and also ${\widehat Q}_\Width$ and
$Q^{\delta}_{\Width}$ (for any $\delta\geq 0$) using corresponding
constructions on $\Xa_i(n)$ and $\Xd_i(n)$. When $n<0$, there are
inclusion maps
$i^+_\Width \colon \Xp_i(n;\Width) \longrightarrow \Xp_i(n)$
and also ${\widehat i}_\Width$ and $i^\delta_{\Width}$ for the other variants.

\begin{lemma}
\label{lemma:TruncateX}
If $\Width$ is sufficiently large, then when $n>0$,
the quotient maps $Q^+_\Width$,
${\widehat Q}_\Width$ and $Q^\delta_\Width$ are quasi-isomorphisms;
when $n<0$, the inclusions $i^+_\Width$, ${\widehat i}_\Width$, and $i^\delta_\Width$
are quasi-isomorphisms.
\end{lemma}

\begin{proof}
This follows quickly from the fact that for all $s$ sufficiently large,
the maps 
$\vertp_s$ (and $\verta_s$, $\vertd_s$)
and $\horp_{-s}$ (and $\hora_{-s}$, $\hord_{-s}$)
are isomorphisms.
\end{proof}

It follows easily from the above that $\Xp_i(n)$ is of $\CFp$-type,
in the sense of Definition~\ref{def:CFpType}.

\subsection{Gradings}
\label{subsec:GradingsOnW}

Let $Y$ be a three-manifold with a null-homologous $K\subset Y$.
Continuing notation from Subsection~\ref{subsec:IntegralSurgeries},
given an integer $N$, let $W'_N(K)\colon Y_N(K) \longrightarrow Y$
denote the induced two-handle cobordism.  Fix also Seifert surface $F$
for $K$, and let ${\widehat F}\subset W'_{N}(K)$ denote the surface
obtained by capping off $F$ in $W'_N(K)$.

In the following statement, note that we continue conventions from
Subsection~\ref{subsec:CFcx}: $\CFd(Y_N(K),[i])$ and $\CFd(Y)$ will
denote complexes which are quasi-isomorphic to the Heegaard Floer
homology complexes for some choice of Heegaard diagram and auxilliary data.

\begin{lemma}
        \label{lemma:GradingsOnW} Suppose that $K$ is a knot in a
        three-manifold $Y$. Then, for any constant $\delta\geq 0$,
        there is another constant $\Width$ with the following
        property. For all sufficiently large $N$, for any fixed
        $[i]\in\Zmod{N}$, there are at most two $\SpinC$ structures on
        $W'_{N}(K)$ whose restriction to $Y_N(K)$ is $[i]$ with the
        property that \begin{equation} \label{eq:SmallestSpinC} \min
        \gr \CFd(Y)-\max\gr \CFd(Y_N,[i]) \leq
        \frac{c_1(\spinc)^2+1}{4}; \end{equation} these are the
        $\SpinC$ structures with $$\langle c_1(\spinc),[{\widehat
        F}]\rangle = 2s\pm N,$$ where $s$ is an integer $\equiv
        i\pmod{N}$, and which satisfy $-N/2\leq s < N/2$.  All other
        $\SpinC$ structures satisfy the inequality \begin{equation}
        \label{eq:BadSpinCStructures} c_1(\spinc)^2 \leq -4N.
        \end{equation} Moreover, if any integral representative $s$
        for $i$ satisfies $|s|>\Width$, then there is a unique
        $\SpinC$ structure satisfying
        Equation~\eqref{eq:SmallestSpinC}, and it is the one for which
        $|\langle c_1(\spinc),[\CapSurf]\rangle|$ is minimal.
\end{lemma}

\begin{proof}
  The set of $\SpinC$ structures over $W'_N(K)$ can be identified with 
  the set of integers via the correspondence which sends $\spinc$ to the 
  integer $s$
  determined by the formula
  $$ \langle c_1(\spinc),[{\widehat F}] \rangle -N = 2s.$$
  Fixing the restriction of $\spinc$ to $Y_N(K)$ corresponds to 
  fixing the equivalence class of $s$ modulo $n$.
  Thus, maximizing  the function $(c_1(\spinc)^2+1)/4$ 
  over the set of $\SpinC$ structures whose restriction to $Y_N(K)$
  corresponds to $i\in\Zmod{N}$ amounts to maximizing the
  quadratic function
  $$
  q(s)=\frac{1}{4}\left(1-\frac{(N+2s)^2}{N}\right)
  $$
  over all integers $s$ with $s\equiv i\pmod{N}$.
  
  This quadratic form is maximized by the representative $i$ for $[i]$
  with $-N<i\leq 0$. The second largest value is smaller by
  $\min(2|i|,2N+2i)$, and all other values are smaller than $q(i)$ by at
  least $\max(2|i|,2N+2i)$.  Thus, if $\spinc$ is any $\SpinC$ structure
  for which the quadratic function takes on neither its maximal or
  next-to-maximal value, we have that
  $$\frac{c_1(\spinc)^2+1}{4}\leq q(i)-\max(2|i|,2N+2i)\leq q(i)-N
  \leq \frac{1}{4}-N.$$
  On the other hand, according to Corollary~\ref{cor:LargeNFloerHomology},
  there is a constant $C_2$ with the property that 
  if we let $\Delta=\min\gr\CFd(Y)-\max\gr\CFd(Y_N,[i])$,
  then
  $$q(i)\leq \Delta+C_2.$$
  Thus, if $N$ is sufficiently large, Inequality~\eqref{eq:SmallestSpinC}
  is violated.

  Moreover, for the final statement, observe that the hypothesis
  that any representative $s$ for $i$ satisfies $|s|>\Width$ is equivalent
  to the assertion that $\min(2|i|,2(N+i))> 2\Width$. Moreover,
  if $\spinc$ is any $\SpinC$ structure for which the quadratic form
  is not minimized,
  $$\frac{c_1(\spinc)^2+1}{4}\leq q(i)-\min(2|i|,2N+2i)
  \leq \Delta + C_2 -\min(2|i|,2N+2i)<\Delta+ C_2-2\Width;$$
  thus the choice of $\Width=C_2/2$ satisfies the final assertion.
\end{proof}

\subsection{The case of $\HFa$.}
\label{subsec:HFa}

Theorem~\ref{thm:IntegerSurgeries} follows from an analysis of
Theorem~\ref{thm:ExactSeq} for a suitable choice of $m$ --
specifically, we choose $m=nk$ where $k$ is a sufficiently large
positive integer.  We shall assume for the time being that $n>0$,
returning to the case of negative surgery coefficients in
Subsection~\ref{subsec:NegativeN}. Our aim in the present subsection is
to prove an analogue of Theorem~\ref{thm:IntegerSurgeries} for $\CFa$,
showing that $\Xa_i(n)$ is quasi-isomorphic to $\CFa(Y_n(K),i)$.
(However, the statement about maps induced by cobordisms will be
relegated to Subsection~\ref{subsec:Maps} below, and the case where
$n<0$ is handled in Subsection~\ref{subsec:NegativeN}.)

Let $\fa_1\colon \CFa(Y_n(K)) \longrightarrow \CFa(Y_{n(k+1)}(K))$ and
$\fa_2\colon \CFa(Y_{n(k+1)}(K))\longrightarrow \CFa(Y;\Z[\Zmod{nk}])$
denote the maps induced by $\fp_1$ and $\fp_2$ from
Equations~\eqref{eq:DefF1} and~\eqref{eq:DefF2} on $\CFa\subset \CFp$,
thought of as the kernel of multiplication by $U$.
We will also fix an identification of chain complexes 
$$\CFa(Y;\Z[\Zmod{nk}])\cong \CFa(Y)\otimes_{\Z}\Z[\Zmod{nk}],$$
and so as to think of $\CFa(Y;\Z[\Zmod{nk}])$ as a direct sum
of complexes
\begin{equation}
\label{eq:DirectSumDecomposition}
\CFa(Y;\Z[\Zmod{nk}])\cong \bigoplus_{s\in\Zmod{nk}} T^s\otimes \CFa(Y).
\end{equation}
(We are using here an identification $\theta$ as in
Equation~\eqref{eq:DefChangeP}, which is well-defined up to an overall
shift by $T^c$ for some constant $c$ which we will fix later. We suppress
the identification $\theta$ to simplify notation.)

Let $W'_{n(k+1)}(K)$ be the natural two-handle cobordism from
$Y_{n(k+1)}(K)$ to $Y$. We abbreviate this cobordism by $W'(k)$.
Recall that a choice of Seifert surface $F$ for $K$ gives rise to a
closed surface ${\widehat F}\subset W'(k)$ with
$$[{\widehat F}]\cm[{\widehat F}]=-n(k+1).$$

Given $s\in\Z$, let $\AYa_s= \CFa(Y_{n(k+1)}(K),s)$; i.e. it is the
summand of $\CFa(Y_{n(k+1)}(K))$ in a $\SpinC$ structure gotten by
restricting a $\SpinC$ structure ${\mathfrak s}$ over $W'(k)$ which
satisfies
$$\langle c_1(\spinc)\cm [{\widehat F}]\rangle -n(k+1)\equiv 2s
\pmod{2n(k+1)}.$$
Moreover, let $\spincx_{s}$ resp. $\spincy_{s}$ denote
the $\SpinC$ structures over the cobordism $$W'(k)\colon
Y_{n(k+1)}(K)
\longrightarrow Y$$ which satisfy
\begin{eqnarray*}
\langle c_1(\spincx_s),[{\widehat F}] \rangle +n(k+1) = 2s
&{\text{and}}&
\langle c_1(\spincy_s),[{\widehat F}] \rangle -n(k+1) = 2s
\end{eqnarray*}
respectively. 

For any $s$ in
\begin{equation}
\label{eq:ChoiceOfs}
-\frac{n(k+1)}{2} \leq s < \frac{n(k+1)}{2},
\end{equation}
then $\spincx_s$ and $\spincy_s$ are the two $\SpinC$ structures over
$W'(k)$ with fixed restriction to $Y_{n(k+1)}(K)$ for which the
function $\frac{c_1(\spinc)^2+1}{4}$ takes on its two largest values
(cf. Lemma~\ref{lemma:GradingsOnW}).

Observe that $\spincx_s+\PD[{\widehat F}] = \spincy_s$, and
hence 
\begin{equation}
\label{eq:spincxMspincy}
{\mathfrak m}(\spincx_{s})-{\mathfrak m}(\spincy_s) = -n(k+1)
\equiv -n \pmod{nk},
\end{equation} and also 
$${\mathfrak
m}(\spincx_{s+n})-{\mathfrak m}(\spincx_s)
\equiv n\pmod{nk}.$$
Thus, if we write
$$\BYa_s= T^s\otimes \CFa(Y) \subset \CFa(Y;\Zmod{m})$$
(with respect to the direct sum decomposition of
Equation~\eqref{eq:DirectSumDecomposition}), then, after multiplying
$\fp_2$ with an overall factor of $T^c$ for some constant $c$, the
components ${\widehat f}_2$ corresponding to the $\SpinC$ structures
$\spincx_s$ and $\spincy_s$ for integers $s$ in the range of
Inequality~\eqref{eq:ChoiceOfs} give maps
\begin{eqnarray}
\label{eq:VertHor}
\vertYa_s\colon \AYa_s \longrightarrow \BYa_s 
&{\text{and}}&
\horYa_s\colon \AYa_s \longrightarrow \BYa_{s+n}
\end{eqnarray}
respectively.

Now, apply Lemma~\ref{lemma:GradingsOnW} in the case where $\delta=0$,
so that $\CFd=\CFa$, and choose $s$ to satisfy
Inequality~\eqref{eq:ChoiceOfs}. In this case, for all sufficiently
large $k$, the lemma combined with the dimension shift formula
(Equation~\eqref{eq:DimensionShift}) proves that $\spincx_s$ and
$\spincy_s$ are the only two $\SpinC$ structures which may induce
non-trivial maps from $\CFa(Y_{n(k+1)}(K),s)$ into $\CFa(Y)$.

In sum, we have identified ${\widehat f}_2$, the map on $\CFa$ induced
from the map $\fp_2$ from Equation~\eqref{eq:DefF2},
with
the map
$$\fa_2'\bigoplus_{s\in \Zmod{n(k+1)}}\AYa_{s} 
\longrightarrow \bigoplus_{s\in\Zmod{nk}} \BYa_s$$
induced by adding all the $\vertYa_s$ and $\horYa_s$
in the range specified by Equation~\eqref{eq:ChoiceOfs}
(compare Equation~\eqref{eq:InvariantF2}).

Moreover,
according to the second statement in Lemma~\ref{lemma:GradingsOnW},
there is an integer $\Width$ (independent of $k$) with the property that
for all $s\geq \Width$, we have that
\begin{eqnarray*}
\vertYa_{-s} \colon \AYa_{-s}\longrightarrow
\BYa_{-s} &{and}&
\horYa_{s} \colon \AYa_{s}\longrightarrow
\BYa_{s+n} 
\end{eqnarray*}
are null-homotopic. This, together with the integer surgeries long
exact sequence (and the fact that $\HFp(S^3_0(K),s)=0$ for $s\geq \Width$ for 
some $\Width$), also shows that
\begin{eqnarray*}
\vertYa_{s} \colon \AYa_{s}\longrightarrow
\BYa_{s} &{and}&
\horYa_{-s} \colon \AYa_{-s}\longrightarrow
\BYa_{-s+n} 
\end{eqnarray*}
are quasi-isomorphisms.
It follows that the mapping cone of $f_2'$
is quasi-isomorphic to the mapping cone of 
$$
\fa_2''\colon \bigoplus_{\{s\in\Z\big| s\equiv i\pmod{n},
  -\Width\leq s\leq \Width\}}
\AYa_{s} \longrightarrow 
\bigoplus_{\{s\in\Z\big| s\equiv i\pmod{n},
  -\Width+n\leq s\leq \Width\}}
\BYa_{s}
$$
obtained by adding all the $\vertYa_s$ and $\horYa_s$
in the given range.

Theorem~\ref{thm:LargeNSurgery} gives
interpretations of these objects in terms of the knot Floer homology.
Indeed, we have identifications (provided that $k$ is sufficiently large)
\begin{eqnarray*}
\begin{CD}
\AYa_{s} @>{\vertYa_{s}}>> \BYa_s \\
@V{\Psi^+_{n(k+1),s}}VV    @VVV    \\
\Aa_{s} @>{\verta_{s}}>> B_s 
\end{CD}
&{\text{and}}&
\begin{CD}
\AYa_{s} @>{\horYa_{s}}>> \BYa_{s+n} \\
@VV{\Psi^+_{n(k+1),s}}V    @VVV    \\
\Aa_{s} @>{\hora_{s}}>> B_{s+n}.
\end{CD}
\end{eqnarray*}
This in turn shows that the mapping cone of $\fa_2''$ is quasi-isomorphic to
the mapping cone of
$$\fa_2'''\colon \bigoplus_{\{s\in\Z\big| s\equiv i\pmod{n},
  -\Width\leq s\leq \Width\}}
\Aa_{s} \longrightarrow 
\bigoplus_{\{s\in\Z\big| s\equiv i\pmod{n},
  -\Width+n\leq s\leq \Width\}}
\Ba_{s}
$$
obtained by adding $\verta_s$ and $\hora_s$ in the given range
(cf. Lemma~\ref{lemma:qIso}).
But this mapping cone is identified 
the truncation $\Xa_i(n,\Width)\subset \Xp_i(n,\Width)$
from Subsection~\ref{subsec:Truncation}.  
Finally, applying Lemma~\ref{lemma:TruncateX},
we see that $\CFa(Y_n,i)$ is quasi-isomorphic to the mapping cone $\Xa_i(n)$.

\subsection{The case of $\CFd$}
\label{subsec:CFd}

A direct application of the argument from
Subsection~\ref{subsec:HFa} proves the following:

\begin{prop}
\label{prop:CFd}
Fix integers $n>0$ and $\delta\geq 0$. Then, there is a constant $\Width$
with the property that $\CFd(Y_n(K),i)$ is quasi-isomorphic to
$\Xd_i(n;\Width)$ (in the notation of Subsection~\ref{subsec:Truncation}).
\end{prop}

\begin{proof}
Apply the proof of Subsection~\ref{subsec:HFa}, only now apply
Lemma~\ref{lemma:GradingsOnW} with $\delta\geq 0$ arbitrary (rather
than $=0$).
\end{proof}

Of course, according to Lemma~\ref{lemma:TruncateX}, the truncation is
unnecessary, and we see that $\HFd(Y_n(K),i)\cong H_*(\Xd_i(n))$.
Moreover, when working with Floer homology with coefficients in a
field, this statement for all $\delta\geq 0$ suffices to prove a
version of Theorem~\ref{thm:IntegerSurgeries} in an ungraded sense.
We do not pursue this direction, but instead turn to gradings to
establish the stronger form of the result stated in
Theorem~\ref{thm:PreciseIntegerSurgeries}.

\subsection{More gradings.}
\label{subsec:MoreGradings}

We turn our attention now to the map $\fp_1$ from
Equation~\eqref{eq:DefF1}, gotten by counting holomorphic triangles in
the four-manifold $X_{\alpha,\gamma,\delta}$ belonging to the Heegaard
triple $(\Sigma,\alphas,\deltas,\betas)$.  This four-manifold has
three boundary components,
$$Y_n(K),\hskip1cm Y_{\gamma,\delta}\cong L(nk,1)\#
(\#^{g-1}(S^2\times S^1)), {\text{~and~}} \hskip1cm Y_{n(k+1)}(K).$$
We abbreviate this four-manifold by $X(k)$.  We will always fix the
canonical $\SpinC$ structure over the boundary component
$Y_{\gamma,\delta}$ (cf. Definition~\ref{eq:CanonicalSpinC}).

For $s\in\Z$, let
$$\Pi^A_s \colon \CFp(Y_{n(k+1)}(K)) \longrightarrow \CFp(Y_{n(k+1)}(K),s)$$
denote the natural projection map.
Consider the map 
$$\fp_1\colon \CFp(Y_n(K))\longrightarrow \CFp(Y_{n(k+1)}(K))=
        \bigoplus_{s\in\Zmod{n(k+1)}}\CFp(Y_{n(k+1)}(K),s),$$
from Equation~\eqref{eq:DefF1}
and let $\fd_1$ denote its restriction to $\CFd(Y_n(K))$.

In this subsection, we prove the following result.

\begin{prop}
\label{prop:RelativelyGraded}
Fix an absolute lift of the relative $\Z$-grading on $\CFp(Y_n(K))$,
and an integer $\delta\geq 0$. There exists a constant $\Width$ so that,
for all sufficiently large $k$, 
there are absolute lifts of the
relative $\Z$-gradings on both 
$\CFd(Y_{n(k+1)}(K),s)\subset \CFd(Y_{n(k+1)}(K))$ 
and
$T^{s}\otimes\CFd(Y)\subset \CFd(Y;\Z[\Zmod{nk}])$ 
for all $|s|\leq \Width$,
with the property that $\Pi^A_s\circ \fd_1$ and also the restriction
of $\fd_2$ to $\CFd(Y_{n(k+1)}(K),s)$ 
have degree zero.
\end{prop}

We give the proof after setting up some terminology and establishing a lemma.

Sometimes, we find it convenient to pass between the absolute $\Q$-gradings
on $\CFd(Y_{n(k+1)}(K))$, $\CFd(Y;\Z[\Zmod{nk}])$ and the absolute gradings
induced from the above proposition (e.g. induced from
the absolute $\Q$-grading on $\CFd(Y_n)$).  We call the 
absolute $\Q$-grading from before the ``old gradings'', and we call
the induced gradings the ``new gradings''.

\begin{lemma}
\label{lemma:GradingsOnX}
Fix a constant $C_0$.  For all
sufficiently large $k$, the following statement holds.  Each
$\SpinC$ structure over $Y_{n(k+1)}(K)$, has at most one
extension $\spinc$ over $X(k)$ whose restriction to
$Y_{\gamma,\delta}$ is the canonical $\SpinC$ structure and for which
\begin{equation}
  \label{eq:MaximalSpinCInX}
  C_0 \leq c_1(\spinc)^2+nk.
\end{equation}
\end{lemma}

\begin{proof}
  It is easy to see that a generator $\Sigma$ for
  $H_2(X(k);\Z)$ has
  $$\Sigma^2=-nk(k+1).$$
  Thus, if $\spinc$ satisfies Inequliaty~\eqref{eq:MaximalSpinCInX},
  then
  $c_1(\spinc)=\alpha\cm
  \PD[\Sigma]$, where $\alpha$ satisfies
  $$\Big|\alpha\Big|\leq \sqrt \frac{nk-C_0}{nk(k+1)} \leq \frac{1}{2}$$
  for all sufficiently large $k$.
  
  Note that any other $\SpinC$ structure which interpolates between
  the same two $\SpinC$ structures on $Y_{n}$ and $Y_{n(k+1)}$ has 
  the form $\spinc+\ell\cm \PD[\Sigma]$ for some integer $\ell\neq 0$; 
  now,
  \begin{eqnarray*}
    (c_1(\spinc+\ell\cm\PD[\Sigma]))^2-c_1(\spinc)^2
    &=&
    4(\ell^2\Sigma\cm\Sigma+\ell\cm \langle c_1(\spinc),[\Sigma]\rangle) \\
    &\leq &
    -4\ell^2 nk (k+1)\left(1-\frac{1}{2\ell}|\alpha|\right) \\
  &\leq& -2 nk(k+1)
  \end{eqnarray*}
  if $k$ is sufficiently large.  In particular, for sufficiently large
  $k$, Inequality~\eqref{eq:MaximalSpinCInX} is violated.
\end{proof}

\vskip.2cm
\noindent{\bf Proof of Proposition~\ref{prop:RelativelyGraded}.}
Fix an absolute lift of the relative $\Z$-grading on $\CFp(Y_n(K))$,
and fix some integer $\delta\geq 0$.
Observe that $\fd_1$ can be decomposed as a sum of homogeneous terms,
indexed by $\SpinC$ structures $\spinc$ over the four-manifold $X(k)$
whose restriction to $Y_{\delta,\gamma}$ is the canonical $\SpinC$ structure.
Each term is homogeneous, and with respect
to the natural $\Q$-gradings on $\CFp(Y_n(K))$ and $\CFp(Y_{n(k+1)}(K))$
(i.e. the ``old gradings''),
they are homogeneous of degree
$$\frac{c_1(\spinc)^2+nk}{4}$$
(compare Equation~\eqref{eq:DimensionShift}; observe that we have a third
boundary component in this four-manifold, where we use a fixed generator
$\Theta_{\delta,\gamma}$).
In view of Corollary~\ref{cor:LargeNFloerHomology} there are constants
$C_1$ and $C_2$ with the property that for all sufficiently large $k$,
\begin{eqnarray*}
    \lefteqn{\min\gr \CFd(Y_{n(k+1)}(K))-\max\gr \CFd(Y_n(K))} \\
    &\geq& 
    \min\gr\CFd(Y_{n(k+1)}(K)) - \max\gr\CFd(Y_{n(k+1)}(K)) \\
    && + \max\gr\CFd(Y_{n(k+1)}(K))-\min\gr\CFd(Y) \\
    && + \min\gr\CFd(Y)-\max\gr\CFd(Y_n(K)) \\
    &\geq& -C_1-C_2+\min_{i\in\Zmod{n(k+1)}} d(n(k+1),i) \\
    && +\min\gr\CFd(Y)-\max\gr\CFd(Y_n(K)) \\
    &\geq & -\frac{1}{4}-C_1-C_2
    +\min\gr\CFd(Y)-\max\gr\CFd(Y_n(K))
  \end{eqnarray*}
Applying Lemma~\ref{lemma:GradingsOnX} with 
$$C_0=-(1+4(C_1+C_2-\min\gr\CFd(Y)+\max\gr\CFd(Y_n)),$$
we see that for all sufficiently large $k$, there is at most
one $\SpinC$ structure which interpolates between given $\SpinC$
structures on $\CFd(Y_n(K))$ and $\CFd(Y_{n(k+1)}(K))$ and which satisfies
the inequality
$$\min\gr \CFd(Y_{n(k+1)}(K))-\max\gr \CFd(Y_n(K)) \leq \frac{c_1(\spinc)^2+nk}{4}.$$
Clearly, if this inequality is not satisfied,
then the $\SpinC$ structure $\spinc$ induces a trivial map
from $\CFd(Y_n(K))$ to $\CFd(Y_{n(k+1)}(K))$.
Thus, in view of Lemma~\ref{lemma:GradingsOnX}, when $k$ is
sufficiently large, the map
$$\Pi^A_s\circ 
\fd_1\colon \CFd(Y_n(K)) \longrightarrow \CFd(Y_{n(k+1)}(K),s)$$
is homogeneous, and hence there is an induced grading on $\CFd(Y_{n(k+1)}(K))$
for which $\fd_1$ is a degree zero map. This induced grading, of course, depends on the intersection form on $X(k)$. 

In the proof of Proposition~\ref{prop:CFd}, the mapping cone of $\fd_2$ 
is identified with the mapping cone of
$$\Dd_{n,i;\Width}\colon \bigoplus_{\{s\in\Z\big| s\equiv i\pmod{n},
  -\Width\leq s\leq \Width\}}
\Ad_{s} \longrightarrow 
\bigoplus_{\{s\in\Z\big| s\equiv i\pmod{n},
  -\Width+n\leq s\leq \Width\}}
\Bd_{s},$$
gotten by adding the maps
\begin{eqnarray*}
\vertd_s\colon \Ad_s \longrightarrow \Bd_s
&{\text{and}}&
\hord_s\colon \Ad_{s}\longrightarrow \Bd_{s+n}
\end{eqnarray*}
in the given range.  Having just endowed $\CFd(Y_{n(k+1)}(K),s)\cong
\Ad_s$ with an absolute grading, there are absolute gradings on
$\Bd_s$ for which $\vertd_s$ have degree zero. We claim that for these
choices of absolute gradings, the map $\hord_s$ is also a
grading-preserving map for $|s|\leq \Width$. This amounts to establishing
that $\vertd_s\circ\Pi^A_s \circ \fd_1$ and
$\hord_{s+n}\circ\Pi^A_{s+n}\circ \fd_1$ are both homogeneous
maps with the same degree (for brevity, we suppress here the
identification of $\Ad_s\cong \CFd(Y_{n(k+1)},s)$ from the notation).
Note first that both these maps are homogeneous, and that their degree
difference is some constant $c$ depending on $i$, $s$, $n$, and $k$
(but independent of the knot $K$). Thus, to verify that $c=0$, it
suffices to verify that $c=0$ in a single model case.  More precisely,
we need a particular $K\subset Y$ and chains $\xi\in\CFp(Y_n(K),i)$ so
that for all sufficiently large $k$,
$\vertp_s\circ\Pi^A_s\circ\fd_1(\xi)$ and
$\horp_{s+n}\circ\Pi^A_{s+n}\circ\fd_1(\xi)$ are homogeneous elements
of the same degree.

To this end, we consider the unknot in $S^3$. Although
$\Pi^A_s\circ\fd_1$ is, in general, difficult to calculate, in the case
of the unknot, it is straightforward. 

Specifically, fix any $\delta\geq 0$ and any 
$\Width \geq 0$.
According to Proposition~\ref{prop:CFd}, for all sufficiently large $k$,
the map $\phi^\delta = (\fd_1,\Hd_1)$ induces
a quasi-isomorphism 
$$\CFd(S^3,i) \longrightarrow \Xd_i(n;\Width).$$
It follows from
Subsection~\eqref{subsec:Unknot} that the induced map on homology must
agree (up to an overall sign) 
with the map $\iota$ from Equation~\eqref{eq:DefIota} or, more
precisely, the restriction of $\iota$ to $\HFd(S^3)$. By taking
$\delta$ sufficiently large, the map $\iota$ has non-trivial
components in both $H_*(\Ad_s)$ and $H_*(\Ad_{s+n})$ Moreover,
provided $s>0$, the map $\vertd_{s}$ induces an isomorphism in
homology, and hence there is some cycle $\xi$ representing a
non-trivial homology class (of arbitrary even degree). Since $\xi$ is
a cycle, so is its image under the quasi-isomorphism with
$\Xd_i(n;\Width)$.  In particular, for any integral $s$, we can choose
$\delta$ sufficiently large for $\Pi^A_s\circ \fd_1(\xi)$ to be
non-trivial, and even for its image under $\vertd_s$ to be a
non-trivial cycle in $\Bd_s$. Since $\fd_1(\xi)$ is a cycle in the
mapping cone, it follows that $\horp_{s-n}\circ \Pi^A_{s-n}\circ
\fd_1(\xi)$ is homologous to $\vertp_s\circ\Pi^A_{s}\circ\fd_1(\xi)$,
and in particular, they are supported in the same degree.  \qed
\vskip.2cm

\subsection{Completion of the proof of Theorem~\ref{thm:PreciseIntegerSurgeries}
  in the case where $n>0$.}
\label{thm:PositiveN}

According to Proposition~\ref{prop:RelativelyGraded}, an
absolute grading on $\CFd(Y_n)$ induces absolute gradings on both
\[\bigoplus_{
-\Width\leq s \leq \Width}\CFd(Y_{n(k+1)},s)\]
and 
\[\bigoplus_{
-\Width+n\leq s \leq \Width} T^s\otimes \CFd(Y) \subset
\CFd(Y;\Z[\Zmod{nk}])
\]
so that the maps $\fd_1$
and $\fd_2$ both are graded maps with degree zero. 
Let 
$$\Pi^B_s\colon \CFd(Y;\Z[\Zmod{nk}])\longrightarrow \CFd(Y)\cong \Bd_s$$
denote projection onto the summand of the form $T^s\otimes \CFd(Y)$.

\begin{lemma}
\label{lemma:GradedHomotopy}
With respect to the above gradings, given $\delta\geq 0$, we have that
for any $k$ sufficiently large and $|s|\leq \Width$, the map $$\Pi^B_s\circ
\Hd_1\colon \CFd(Y_n) \longrightarrow \Bd_s\cong
\CFd(Y)$$ is a homogeneous map of degree $+1$.
\end{lemma}

\begin{proof}
$\Hd_1$ splits as a sum of homotopy classes of maps of holomorphic
quadrilaterals. These homotopy classes in turn correspond to $\SpinC$
structures over the cobordism from $Y_n$ to $Y$. Since the
intersection form on the composite cobordism $X(k)\cup_{Y_{n(k+1)}}
W'(k)$ is negative-definite, it follows that, $\Hd_1$ is a sum of
homogeneous maps, whose top order part corresponds to $\SpinC$
structures over the composite cobordism whose first Chern classes have
maximal square. 

If the $\spinc$ component of $\Hd_1$ is non-trivial,
then we claim that $\spinc|X(k)$ is a maximal $\SpinC$ structure over
$X(k)$ in the sense that it appears in a component of $\fd_1$; and
also, $\spinc|W'(k)$ is one of the top two $\SpinC$ structures over
$W'(k)$ in the sense that it appears as a component of $\fd_2$.

To this end, note that with respect to the usual gradings on $\CFp(Y_n)$ and
$\CFp(Y;\Z[\Zmod{nk}])$, the component of $\Hd_1$ using the $\SpinC$
structure $\spinc$ is a homogeneous map of degree
$$1 + \frac{c_1(\spinc|X(k))^2+nk}{4} + \frac{c_1(\spinc|W'(k))^2+1}{4}.$$
Thus, this map is non-trivial only if
\begin{equation}
  \label{eq:NonTrivialMap}
  \min\gr\CFp(Y)-\max\gr\CFp(Y_n)\leq
  1 + \frac{c_1(\spinc|X(k))^2+nk}{4} + \frac{c_1(\spinc|W'(k))^2+1}{4} 
\end{equation}
However, 
$$  1 + \frac{c_1(\spinc|X(k))^2+nk}{4} + \frac{c_1(\spinc|W'(k))^2+1}{4} 
\leq 
1 + \frac{c_1(\spinc|X(k))^2+nk}{4} + \frac{1}{4},
$$
and since the left-hand-side of Equation~\eqref{eq:NonTrivialMap}
is independent of $k$,
it follows from Lemma~\ref{lemma:GradingsOnX} that
(provided that $k$ is sufficiently large)
$\spinc|X(k)$ is one of the distinguished $\SpinC$ structures on
$X(k)$ satisfying Inequality~\eqref{eq:MaximalSpinCInX}.

Similarly, 
\begin{eqnarray*}
1 + \frac{c_1(\spinc|X(k))^2+nk}{4} + \frac{c_1(\spinc|W'(k))^2+1}{4} 
&\leq &
1 + \frac{nk}{4} + \frac{c_1(\spinc|W'(k))^2+1}{4}, \\
\end{eqnarray*}
For $k$ sufficiently large, this means also that $\spinc|W'(k)$ is one
of the $\SpinC$ structures satisfying
Equation~\eqref{eq:SmallestSpinC}; for if it were not, then by
Equation~\eqref{eq:BadSpinCStructures},
$$
1 + \frac{nk}{4} + \frac{c_1(\spinc|W'(k))^2+1}{4} 
\leq 
\frac{5-4n}{4} -\frac{3nk}{4},
$$
which, if $k$ is sufficiently large, violates Equation~\eqref{eq:NonTrivialMap}.

We have thus verified that if the $\spinc$-component of $\Hd_1$ is
non-trivial, then the restriction of $\spinc$ to $X(k)$ and
$W(n(k+1))$ are the ones which induce the map $\Pi_{s}^A\circ \fd_1$ and
$\fd_2|_{\Ad_{s}}$ (with $|s|\leq \Width$) respectively.

Thus, it follows that $\Hd_1$ is a homogeneous map from $\CFp(Y_n(K))$
to $\CFp(Y;\Z[\Zmod{nk}])$ with respect to the ``new gradings''
(from Proposition~\ref{prop:RelativelyGraded}) with degree
$+1$.
\end{proof}

\vskip.2cm
\noindent{\bf{Proof of Theorem~\ref{thm:PreciseIntegerSurgeries} for $n>0$.}}
We endow the mapping cone $\MCone(\Dd_n)$ with an absolute grading
where $\BigAd$ has the inherited absolute grading, while $\BigBd$ has
the inherited grading, shifted up by one (this is done so that the
differential on the mapping cone is a graded map which shifts grading
down by $-1$).  

According to Proposition~\ref{prop:CFd}, the chain complex
$\CFd(Y_n,i)$ is quasi-isomorphic to the mapping cone
$\Xd_{i}(n;\Width)=\MCone(\Dd_{n,i;\Width})$ via a 
quasi-isomorphism 
induced by
$\phi^{\delta}=(\fd_1,\Hd_1)$ which, according to the above lemma, is
a grading-preserving map.

In view of Lemma~\ref{lemma:TruncateX}, this gives, for any
$\delta\geq 0$ a grading-preserving isomorphism of $\HFd(Y_n,i)$ with
$H_*(\Xd_i(n))$.  It follows from Lemma~\ref{lemma:HFdLemma}
that in fact $\HFp(Y_n,i)\cong H_*(\Xp_i(n))$.
\qed

\subsection{Negative integral surgeries}
\label{subsec:NegativeN}

In fact, the above proof of Theorem~\ref{thm:PreciseIntegerSurgeries}
for positive surgery coefficients 
adapts with minor changes to the case of negative surgery coefficients.
We outline the changes presently.

\vskip.2cm
\noindent{\bf{Proof of Theorem~\ref{thm:PreciseIntegerSurgeries} for negative
surgery coefficients.}}
We will continue to assume $n>0$, only now, we shall consider
the three-manifold $Y_{-n}(K)$.

We start by considering the changes to Subsection~\ref{subsec:HFa}.
In the present case, we apply Theorem~\ref{thm:ExactSeq}, only now 
focusing on the maps
$$\begin{CD}
\CFp(Y_{n(k-1)}) @>{\fp_2}>>
\CFp(Y;\Z[\Zmod{nk}])
@>{\fp_{3}}>>
\CFp(Y_{-n}(K)).
\end{CD}$$ 
Where again we take $k$ to be sufficiently large.

In place of Equation~\eqref{eq:VertHor}, we have
$$
{\mathfrak m}(\spincx_{s})-{\mathfrak m}(\spincy_s) = -n(k-1)
\equiv n\pmod{nk}.
$$
Thus, it follows that, unlike the conventions specified in
Equation~\eqref{eq:VertHor}, the components of $\fp_1$
corresponding to $\spincx_s$ and $\spincy_s$ are maps
\begin{eqnarray*}
\vertYa_s\colon \AYa_s \longrightarrow \BYa_s 
&{\text{and}}&
\horYa_s\colon \AYa_s \longrightarrow \BYa_{s-n}.
\end{eqnarray*}
Proceeding as in Section~\ref{subsec:HFa},
for all sufficiently large $k$, we identify the mapping cone of $\fp_2$ with 
the mapping cone of
$$\fa_2''\colon \bigoplus_{\{s\in\Z\big| s\equiv i\pmod{n},  -\Width\leq s 
\leq \Width-n\}}
\AYa_{s} \longrightarrow 
\bigoplus_{\{s\in\Z\big| s\equiv i\pmod{n},  \Width-n\leq s \leq \Width-n\}}
\BYa_{s} 
$$
This is then identified with the mapping cone $\Xa_i(-n;\Width)$,
via Theorem~\ref{thm:LargeNSurgery} as before.

Indeed, the modifications from Subsection~\ref{subsec:CFd} carry over,
to prove the analogue of Proposition~\ref{prop:CFd}, using
surgery coefficient $-n$.

When studying gradings (cf. Subsection~\ref{subsec:MoreGradings}),
now, we study the map $\fp_3$ belonging to the cobordism $W_{-n}(K)$,
rather than the map $\fp_1$, since now $W_{-n}(K)$ is a
negative-definite cobordism.

Of course, for fixed $\delta\geq 0$ and $k$ large enough, $\fd_2$
remains a homogeneous map according to Lemma~\ref{lemma:GradingsOnW},
as before.  Moreover, we claim that the same holds for $\fd_3$. 
To see
this, consider the cobordism $W_{-n}(K)\colon Y\longrightarrow
Y_{-n}(K)$.  Now, (cf.  Equation~\eqref{eq:InvariantF3})
$\fd_3(\xi\otimes T^{i})$ consists of a sum of maps associated to
$\SpinC$ structures $\spinc$ which differ by addition of $k\cm
[\CapSurf]$ (where here $\CapSurf$ represents a generator for $H_2(W_{-n}(K))$,
and in particular $[\CapSurf]\cm [\CapSurf]=-n$):

\begin{lemma}
  \label{lemma:DistinguishThree}
  Fix an integer $\delta \geq 0$ and a constant $C_0$.
  For all sufficiently large $k$, the following holds.
  For each $\spinc_0\in\SpinC(W_{-n}(K))$, there is 
  at most one $\SpinC$ structure $\spinc\in\spinc_0+k\PD[\CapSurf]\Z$
  for which 
  $$C_0\leq c_1(\spinc)^2+nk.$$
\end{lemma}

\begin{proof}
  Let $\spinc$ be a $\SpinC$ structure as above. In this case,
  $c_1(\spinc)=\alpha\cm\PD[\CapSurf]$, where $\alpha$ satisfies
  $$|\alpha|\leq \sqrt{\frac{nk-C_0}{n}} \leq 2\sqrt{k} $$
  if
  $k$ is sufficiently large. But then, for any other $\SpinC$
  structure, of the form $\spinc_0+(k\ell)\PD[\CapSurf]$ (with $\ell\in\Z$),
  we have that 
  \begin{eqnarray*}c_1(\spinc')^2-c_1(\spinc)^2 
    &\leq &
    4(k\langle c_1(\spinc),[\CapSurf]\rangle -k^2 n)\\
    &\leq &
    -4k^2n (1-\frac{2}{\sqrt k}) \\
    &\leq&
    -nk^2,
  \end{eqnarray*}
        if $k$ is sufficiently large.
\end{proof}

Now, if $$\fd_{W_{-n}(K),\spinc}\colon \CFd(Y) \longrightarrow
\CFd(Y_{-n}(K))$$ is non-trivial, then according to
Equation~\eqref{eq:DimensionShift}, $$c_1(\spinc)^2\geq -1 +
4(\min\gr\CFd(Y_{-n}(K))-\max\gr\CFd(Y)),$$ so the above lemma shows
that for large enough $k$, there is a unique $\SpinC$ structure which
contributes non-trivially to $\fd_3(\cdot\otimes T^i)$.  (This is the
analogue of Lemma~\ref{lemma:GradingsOnX}.)

Thus, a grading on $\CFp(Y_{-n}(K))$ induces gradings on
$\CFp(Y;\Z[\Zmod{nk}])$ and from there also on $\CFp(Y_{n(k-1)})$ (note
that $\fd_2$ remains a graded map for all sufficiently large $k$ by
Lemma~\ref{lemma:GradingsOnW}).  To check that these gradings are
compatible, now, we use the example of the unknot with framing $-n$.

Finally, it remains to show that quasi-isomorphism
$$\psi^\delta\colon \MCone(\fd_2) \longrightarrow \CFp(Y_{-n}(K))$$
defined by $\psi^\delta(x,y)=\Hd_{2}(x)+\fd_3(y)$ is a relatively
graded map. This follows from an analogue of
Lemma~\ref{lemma:GradedHomotopy}. 

To this end, note that $\Hd_2\colon \CFd(Y_{n(k-1)})\longrightarrow
\CFd(Y_{-n})$ is defined as a count of quadrilaterals.
Thus, if the map
$$\Hd_2\colon \CFd(Y_{n(k-1)}(K)) \longrightarrow \CFd(Y_{-n}(K))$$
is non-trivial, then
$$\min\gr\CFp(Y_{-n}(K)) - \max\gr\CFp(Y_{n(k-1)}(K))\leq
1+\frac{c_1(\spinc|W))^2+c_1(\spinc|W''(k))^2+2}{4},$$
where here
$W''(k)$ indicates the cobordism $W'_{n(k-1)}\colon
Y_{n(k-1)}\longrightarrow Y$.  But the left-hand-side is bounded below
by constants (independent of $k$) plus $-nk/4$ (cf.
Corollary~\ref{cor:LargeNFloerHomology}).  It follows now from
Lemma~\ref{lemma:DistinguishThree} that $\spinc|W''(k)$ is one of the
distinguished $\SpinC$ structures satisfying the inequality in
Lemma~\ref{lemma:DistinguishThree}; it also follows from
Lemma~\ref{lemma:GradingsOnW} that $\spinc|W(k)$ is one of the
distinguished $\SpinC$ structures satisfying
Equation~\eqref{eq:SmallestSpinC}. Thus, $\Hd_2$ shifts grading up by $-1$,
and $\psi^\delta$ is a relatively graded map.

\subsection{The case where $n=0$}
\label{subsec:ZeroSurgery}

We have not treated the case where $n=0$, as it follows more quickly
from the existing exact sequences. Specifically, the integer surgeries
long exact sequence (Theorem~\ref{HolDiskTwo:thm:ExactP}
of~\cite{HolDiskTwo}) shows that $\HFp(Y_0(K),i)$ is identified with
the homology of the mapping cone of
$$\vertp_i+\horp_i\colon \Ap_i\longrightarrow \Bp.$$
Indeed,
considering the integer surgeries exact sequence with twisted
coefficients (Theorem~\ref{HolDiskTwo:thm:ExactPTwist}
of~\cite{HolDiskTwo}), we see that $H_*(\Xp_i(0))$ calculates
$\uHFp(Y_0(K),i)$.

\subsection{Maps induced by cobordisms}
\label{subsec:Maps}

The mapping cone $\Xp(n)$ provides also a model for the maps induced
by the cobordism from $Y$ to $Y_n(K)$, according to
Theorem~\ref{thm:Maps}. 
The proof of this result follows quickly from our
proof of Theorem~\ref{thm:PreciseIntegerSurgeries}.

\vskip.2cm
\noindent{\bf{Proof of Theorem~\ref{thm:Maps}.}}
Observe that in Theorem~\ref{thm:ExactSeq} (in the
notation of Section~\ref{sec:ExactSeq}), the following square
commutes:
$$
\begin{CD}
\HFp(Y;\Z[\Zmod{m}]) @>{\Fp{3}}>> \HFp(Y_n(K)) \\
@VV{=}V @V{\Phi^+}VV \\
\HFp(Y;\Z[\Zmod{m}])@>{I}>> H_*(\MCone(\fp_2)),
\end{CD}
$$ where here $\Phi^+$ is the map in homology induced by the
quasi-isomorphism $\phi^+$ of Equation~\eqref{eq:DefPhi}, $I$ is the
map on homology induced by the canonical map $\iota$ to the mapping
cone, and $\Fp{3}$ is the map on homology induced by the map $\fp_3$
from Equation~\eqref{eq:DefF3}.

Let $m=|n|k$ with $k$ sufficiently large, as we have done above.  Fix
also $s\in\Z$ and $\delta\geq 0$. Although the restriction of
$$\fd_3\colon \CFd(Y;\Z[\Zmod{nk}])\longrightarrow \CFd(Y_n(K))$$
to
$T^s\otimes\CFd(Y)\subset \CFd(Y;\Z[\Zmod{nk}])$ in principle is a sum of
maps induced by $\SpinC$ structures over $W_{n}(K)$; all these are in
the $k\PD[\CapSurf]\cm \Z$-orbit of a given one (cf.
Equation~\eqref{eq:InvariantF3}). Thus, for fixed $\delta$ and $k$
sufficiently large, it is easy to see that there is in each such
orbit, a single chain map which could conceivably be non-trivial, and
that is the map $\Fp{W,\spinc}$ induced by
the $\SpinC$ structure with $\langle c_1(\spinc),[\CapSurf]\rangle + n
=2s$.

Thus, we have the following commutative square:
$$
\begin{CD}
  \HFd(Y) @>{F^{\delta}_{W,\spinc}}>> \HFd(Y_n(K)) \\
  @V{=}VV @V{\Phi^\delta}VV \\
  T^s\otimes \HFd(Y) @>{I}>> H_*(\Xd(n;\Width))
\end{CD}
$$
(here, $\Width$ is a truncation as usual, which now we choose so that
$\Width>s$).  Note that $\delta$ here is arbitrary, and all maps are
graded.  Thus, Lemma~\ref{lemma:HFdLemma} provides the canonical
extensions making the following square commute:
$$
\begin{CD}
  \HFp(Y) @>{\Fp{W,s}}>> \HFp(Y_n(K)) \\
  @V{=}VV @V{\cong}VV \\
  T^s\otimes \HFp(Y) @>{I}>> H_*(\Xp(n)),
\end{CD}
$$
proving the theorem.
\qed
\vskip.2cm

\subsection{Generalizations}
\label{subsec:Generalize}

The proof of Theorem~\ref{thm:IntegerSurgeries} applies immediately in
any context where $Y$ has a relatively $\Z$-graded Floer homology.

Specifically, let $Y$ be a closed, oriented three-manifold, and fix a
$\SpinC$ structure $\spinct$ over $Y$ with torsion first Chern class.
Let $K\subset Y$ be a null-homologous knot in $Y$, and fix a Seifert
surface $F\subset Y$ for $K$. 

In this case, the space of $\SpinC$ structures over $Y_n(K)$ which are
$\SpinC$-cobordant to $\spinct$ over the cobordism $W'_n(K)$ is
identified with $\Zmod{n}$, as in Lemma~\ref{lemma:IdentifySpinC}. For $i\in\Zmod{n}$, let
$\CFp(Y_n(K),i,\spinct)$ denote the corresponding summand of
$\CFp(Y_n(K))$. In this case, we let $\Ap_{s,\spinct}$ denote the
corresponding chain complex associated to the knot filtration on
$\CFp(Y,\spinct)$, and let $\Bp_{s,\spinct}$ denote $\CFp(Y,\spinct)$.
Write
\begin{eqnarray*}
\BigAp_{i,\spinct}=\bigoplus_{\{s\in\Z\big|s\equiv i\pmod{n}\}}\Ap_{s,\spinct} &{\text{and}}&
\BigBp_{i,\spinct}=\bigoplus_{\{s\in\Z\big|s\equiv i\pmod{n}\}}\Bp_{s,\spinct} 
\end{eqnarray*}

\begin{theorem}
\label{thm:Generalize}
Fix a $\SpinC$ structure $\spinct$ over $Y$ whose first Chern class is
torsion, and a null-homologous knot $K\subset Y$.  For each $i\in\Zmod{n}$,
the mapping cone
$\Xp_{i,\spinct}(n)$ of 
$$\Dp_{n,i,\spinct}\colon \BigAp_{i,\spinct} \longrightarrow
\BigBp_{i,\spinct},$$
is identified with $\CFp(Y_n(K),i,\spinct)$
as graded $\Z[U]$-modules.
\end{theorem}

Of course, when $n=\pm 1$, $\SpinC(Y_n(K))\cong \SpinC(Y)$, and there
is no additional choice of $i\in\Zmod{n}$. In this case, we write
simply $\Xp_\spinct(n)$ for the mapping cone.

The theorem follows at once from the methods of this section:
specifically, our proof of Theorem~\ref{thm:PreciseIntegerSurgeries}
uses the fact that $Y$ is an integer homology three-sphere only in
that $\CFp(Y)$ has a relative $\Z$-grading. 

\section{Sample Calculations}
\label{sec:Calculations}

We give some calculations to illustrate the techniques of this
paper. In Subsection~\ref{subsec:T34}, we use
Theorem~\ref{thm:PreciseIntegerSurgeries} to calculate the Heegaard
Floer homology groups of $\pm 1$-surgeries on the torus knot $T_{3,4}$
in the three-sphere. (Of course, $\pm 1$ surgeries on torus knots are
Brieskorn spheres, and as such their Heegaard Floer homology groups
can be calculated using the techniques of~\cite{SomePlumbs}; however,
we find it instructive to sketch here the calculations using present
methods.)  In Subsection~\ref{subsec:EulerOne}, we include a
calculation of the Heegaard Floer homology of a non-trivial circle
bundle over a closed, oriented two-manifold. 

\subsection{Surgeries on the torus knot $T_{3,4}$}
\label{subsec:T34}

Sufficiently large surgeries on torus knots have particularly simple
Heegaard Floer homology groups -- in each $\SpinC$ structure,
$\HFp(S^3_N(K))\cong \InjMod{}$, i.e. these spaces are ``Heegaard
Floer homology lens spaces'', or $L$-spaces in the terminology
from~\cite{NoteLens}. As such, the filtered homotopy type of the
knot filtration of torus knots is immediately determined by the
Alexander polynomial, cf. Theorem~\ref{NoteLens:thm:FloerHomology}
of~\cite{NoteLens}. 

\begin{figure}
\mbox{\vbox{\epsfbox{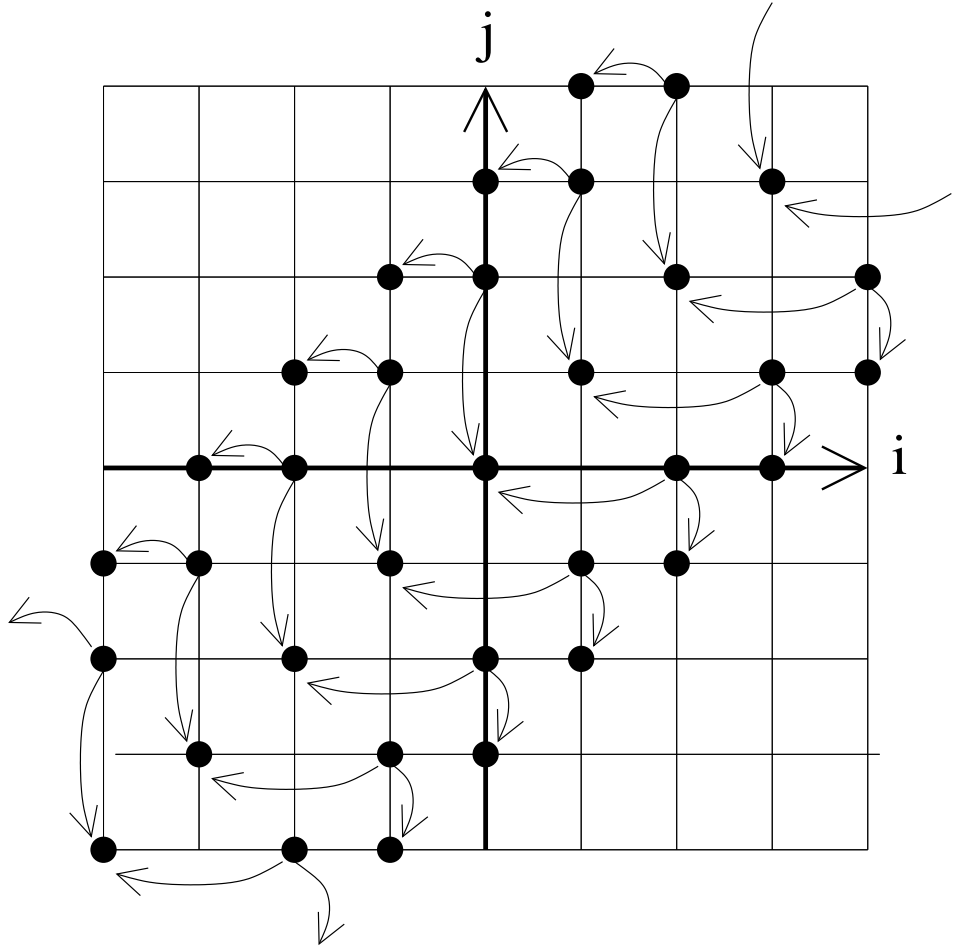}}}
\caption{\label{fig:T34}
{\bf The doubly-filtered knot complex for the torus knot $T_{3,4}$.} We have
illustrated here the plane, representing filtration levels of
generators.  A dot at a lattice point $(i,j)$ represents a group $\Z$
with filtration $(i,j)$. The five generators of $C$ as a $\Z[U,U^{-1}]$-module
are represented by the five dots on the vertical axis.  }
\end{figure}

For example, that result shows that the filtered 
chain complex $C$ for the torus knot $K=T_{3,4}$ has five generators
as a $\Z[U,U^{-1}]$-module,
$\{x_i\}_{i=1}^5$, with 
\begin{eqnarray*}
\Filt(x_1)=(0,3) && \gr(x_1)=0\\
\Filt(x_2)= (0,2) && \gr(x_2)=-1\\
\Filt(x_3)= (0,0) && \gr(x_3)=-2 \\
\Filt(x_4)= (0,-2) && \gr(x_4)=-5\\
\Filt(x_5)= (0,-3) && \gr(x_5)=-6.
\end{eqnarray*}
Moreover, the non-trivial differentials are
\begin{eqnarray*}
\partial x_2 = U\cm x_1+ x_3 &{\text{and}}&
\partial x_4 = U^2\cm x_3 + x_5.
\end{eqnarray*}
The filtrations and the differentials are illustrated in Figure~\ref{fig:T34}.

From this, we quickly see the following:

\begin{lemma}
  \label{lemma:MapsT34} $H_*(\Ap_s)\cong H_*(\Bp_s)\cong \InjMod{}$.
  Moreover, under this identification, 
  the maps on homology induced by
  \begin{eqnarray*}
    \horp_{s}\colon \Ap_{s}\longrightarrow \Bp_{s+n}
    &\text{and}&
    \vertp_{-s}\colon \Ap_{-s}\longrightarrow \Bp_{-s}
  \end{eqnarray*} 
  are identified with multiplication by
$$\left\{\begin{array}{ll}
    U^3 & {\text{if $s=2$}} \\
    U^2 & {\text{if $s=1$}} \\
    U & {\text{if $s=0,-1,-2$}} \\
    1 & {\text{if $s<-2$}}
    \end{array}\right.$$
\end{lemma}

\begin{proof}
  From the above description of $C$, we see that $\Ap_s=\Bp_s$ and
  $\Ap_{-s}=\Bp_{-s-1}$ whenever $s>2$. It remains to consider the
  remaining cases with $|s|\leq 2$.  For example, $\Ap_{1}$ contains
  three additional generators not contained in $\Bp_1$: $U x_1$, $U^2
  x_1$, and $U x_2$. However, these latter two generators do not
  contribute to $H_*(\Ap_1)$: $\partial U x_2=U^2 x_1$ (in $\Ap_1$).
  The bottom-most non-trivial homology class in $H_*(\Ap_1)$ is
  represented by any of the three cycles $U x_1$, $x_3$, or $U^{-2}
  x_5$. However, these cycles are boundaries in $\Bp_1$. It follows
  quickly that $\vertp_1\colon \InjMod{}\cong H_*(\Ap_1)
  \longrightarrow H_*(\Bp_1)\cong \InjMod{}$ is modeled on
  multiplication by $U$. 

  The other cases follow similarly.
\end{proof}

\begin{figure}
\mbox{\vbox{\epsfbox{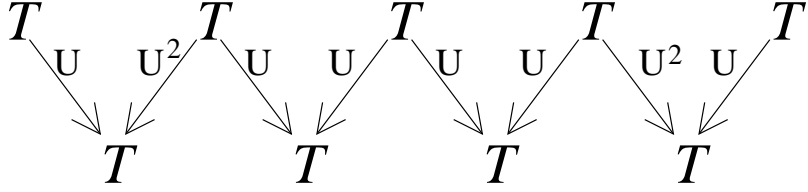}}}
\caption{\label{fig:T34P1}
{\bf Chain complex for $S^3_{+1}(T_{3,4}).$} The results of
Lemma~\ref{lemma:MapsT34} can be summarized by this diagram, which
represents the chain complex $\Xp(1)$. (Note that we have illustrated
here a truncated portion; which carries all of the homology of $\Xp(1)$.)
The arrows represent the differential;
the symbol ${\it T}$ represents the module $\InjMod{}$.
This diagram contains information about gradings as well:
the relative grading is characterized by the requirement that the differential
drops grading by one.}
\end{figure}

The chain complex for $\Xp(1)$ is summarized in Figure~\ref{fig:T34P1}.

In view of Lemma~\ref{lemma:MapsT34}, the map on homology 
$H_*(\BigAp) \longrightarrow H_*(\BigBp)$ induced by $\Dp$ is surjective,
and hence $H_*(\Xp(1))$ can be identified with its kernel.
Indeed,  $H_*(\Xp(1))$
consists of a submodule isomorphic to $\Ap_0$, and also four
additional generators represented by $1\in\InjMod{}\cong H_*(\Ap_s)$
for $0<|s|\leq 2$. According to the grading conventions, the two
bottom-most generators of $\Ap_{\pm 1}$ are supported in degree $-2$
and the two in $\Ap_{\pm 2}$ are supported in degree $0$. Thus, in
light of Theorem~\ref{thm:PreciseIntegerSurgeries}, these calculations
show that:
$$\HFp(S^3_{1}(T_{3,4})) \cong \InjMod{-2}\oplus \Z^2_{(-2)} \oplus
\Z^2_{(0)}.$$
Here, $\InjMod{d}$ denotes a copy of $\InjMod{}$,
thought of as a graded $\Z[U]$-module, graded so that its element of
lowest degree is supported in degree $d$; while $\Z_{(d)}$ denotes
$\Z\cong \Z[U]/U$, supported in degree $d$.

\begin{figure}
\mbox{\vbox{\epsfbox{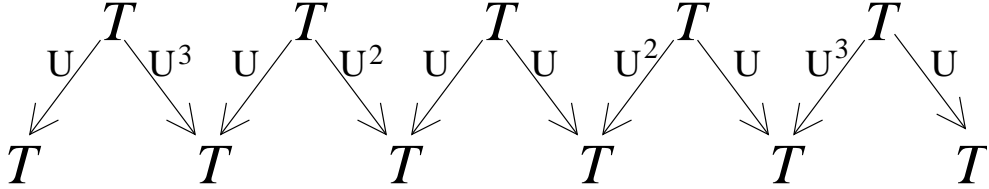}}}
\caption{\label{fig:T34N1}
{\bf Chain complex for $S^3_{-1}(T_{3,4}).$} This is an illustration
of the complex $\Xp(-1)$, with the same conventions
as in Figure~\ref{fig:T34P1}.}
\end{figure}

For $(-1)$-surgery on $T_{3,4}$, 
Lemma~\ref{lemma:MapsT34} gives rise to the chain
complex for $\Xp(-1)$ illustrated in Figure~\ref{fig:T34N1}.

In this case, we can calculate:
$$\HFp(S^3_{-1}(T_{3,4})) \cong \InjMod{0}\oplus \Z_{(-1)} \oplus
\Z^2_{(-3)}\oplus \Z^2_{(-7)}.$$
(These can be compared with the
calculations from~\cite{SomePlumbs}, bearing in mind that
$S^3_{+1}(T_{3,4})\cong -\Sigma(3,4,11)$ and $S^3_{-1}(T_{3,4})\cong
\Sigma(3,4,13)$, where $\Sigma(p,q,r)$ denotes the Brieskorn sphere
with multiplicities $p$, $q$, and $r$, given its orientation as the
boundary of a negative-definite plumbing.)

%
%
%

\subsection{Non-trivial circle bundles over Riemann surfaces}
\label{subsec:EulerOne}

Recall that for any closed, oriented three-manifold, the groups
$U^d\cm \HFp(Y)\subset \HFp(Y)$ stabilize for sufficiently large $d$,
and hence we can define the ``reduced Floer homology group''
$\HFpRed(Y)=\HFp(Y)/U^d\cm \HFp(Y)$ for all sufficiently large $d$.
As an illustration of Theorem~\ref{thm:IntegerSurgeries}, or more
precisely, its generalization in Theorem~\ref{thm:Generalize}, we
calculate the reduced Floer homology of a non-trivial circle bundle
over a Riemann surface of genus $g$, with coefficients in $\Zmod{2}$.

To explain how this is done, recall that there is a genus $g$ fibered
knot in $\#^{2g}(S^2\times S^1)$, the ``Borromean knot'', constructed
as follows. Consider the Borromean link, with a distinguished
component.  Form the $g$-fold connected sum of the Borromean link,
where the connected sum is performed along the distinguished
components. Then, form zero-surgeries on all the remaining components,
viewing the distinguished component as a knot inside
$\#^{2g}(S^2\times S^1)$.  The three-manifold which is an Euler number
$n$ circle bundle over the genus $g$ Riemann surface is obtained from
$n$-framed surgery on this knot.

The knot Floer complex for the Borromean knot is calculated in
Proposition~\ref{Knots:prop:KnotHomology} of~\cite{Knots}.  It
is proved there that $C\cong \Wedge^* H^1(\Sigma;\Z)\otimes
\Z[U,U^{-1}]$. The $\Z\oplus\Z$ filtration is given by
\begin{equation}
  \label{eq:KnotHomologyBorrKnot}
  C\{i,j\}=U^{-i}\otimes \Wedge^{g-i+j}H^1(\Sigma;\Z),
\end{equation}
and the entire group $C\{i,j\}$ is supported in dimension $i+j$.
All the differentials on this knot complex are trivial, and indeed,
so are all the higher differentials; i.e. 
$$\HFinf(\#^{2g}(S^2\times S^1))\cong \bigoplus_{i,j} C\{i,j\}.$$

Indeed, these groups have some additional structure which can be
identified.  Specifically, recall (cf.
Subsection~\ref{HolDisk:subsec:DefAct} of~\cite{HolDisk}) that if $Y$
is an oriented three-manifold, then there is an action of $\Wedge^*
H_1(Y;\Z)/\Tors$ on all the versions of Floer homology. For the case
where $Y\cong \#^{2g}(S^2\times S^1)$, under the identification
$H_1(\Sigma;\Z)\cong H_1(\#^{2g}(S^2\times S^1))$ and
$\HFinf(\#^{2g}(S^2\times S^1))\cong \Wedge^* H_1(\Sigma;\Z)\otimes_\Z
\Z[U,U^{-1}]$ described above, the action of $\gamma\in
H_1(\Sigma;\Z)$ is given by the formula
\begin{equation}
\label{eq:HOneAction}
\gamma\cm (\omega\otimes U^j)=(\iota_\gamma\omega)\otimes U^j +
\PD(\gamma)\wedge \omega \otimes U^{j+1},
\end{equation}
where here $\iota_\gamma$ denotes contraction (cf. the proof of
Theorem~\ref{Knots:thm:SigmaTimesSOne} of~\cite{Knots}).

It is significantly easier to work with coefficients in $\Field=\Zmod{2}$,
because in that case, the induced homotopy equivalence between
$C\{i\geq 0\}$ and $C\{j\geq 0\}$ takes a particularly simple form.
Write $\rC$ for $C\otimes\Field$ (and writing, 
for example, $\rC\{i,j\}$ for $C\{i,j\}\otimes \Field$.
In this case, we have the following:

\begin{prop}
\label{prop:HomEq}
Consider the natural homotopy equivalence 
$$\phi\colon C\{i\geq 0\}\longrightarrow C\{j\geq
0\}.$$ 
Its induced map ${\overline \phi}\colon \rC\{i\geq 0\}
\longrightarrow \rC\{j\geq 0\}$ sends $\rC\{i,j\}$ to 
$\rC\{j,i\}$.
\end{prop}

The above proposition does not hold over $\Z$. A description of the
involution over $\Z$ is given in~\cite{JabukaMark}.

In fact we will be able to calculate the homotopy equivalence 
$\rPhi$, but first, we need a lemma.

\begin{lemma}
\label{lemma:AutHF}
Consider the graded
module 
$$M=\Wedge^*H^1(\Sigma;\Z)\otimes \Z[U,U^{-1}]/\Z[U]$$
over the ring $\Wedge^* H_1(\Sigma;\Z)\otimes_\Z \Z[U]$
with the action defined in Equation~\eqref{eq:HOneAction}.
The only automorphisms of $M$ are multiplication by $\pm 1$.
Similarly, if we consider ${\overline M}=M\otimes\Field$ 
over the ring $\Wedge^* H_1(\Sigma;\Z)\otimes_\Z \Field[U]$,
where $\Field=\Zmod{2}$, the only automorphism of ${\overline M}$
is the identity map.
\end{lemma}

\begin{proof}
  Let $\phi$ be an automorphism of $M$ as a module over $\Wedge^*
  H_1(\Sigma;\Z)\otimes_\Z \Z[U]$.  Consider the filtration of $M$ by
  the submodules $\{M_i\}_{i=1}^\infty$ where $M_i$ is the set of
  elements annihilated by $U^i$.  Clearly, these submodules exhaust
  $M$, and $\phi_i$ preserves this filtration.
  
  Consider the summand $\Z\cong \Wedge^{2g} H^1(\Sigma;\Z)\otimes
  U^{-i} \subset M_i$.  This consists of those elements of $M_i$ with
  maximal grading, and hence it is preserved by $\phi$. Thus, $\phi$
  must induce multiplication by $\pm 1$ on this summand.  However,
  this summand generates all of $M_i$ as a module over $\Wedge^*
  H_1(\Sigma;\Z)\otimes \Z[U]$. It follows at once that $\phi$ induces
  multiplication by $\pm 1$ on all of $M_i$.

  The same argument shows that an automorphism ${\overline \phi}$
  of ${\overline M}$ is the identity map.
\end{proof}

\vskip.2cm
\noindent{\bf{Proof of Proposition~\ref{prop:HomEq}.}}
We describe the homotopy equivalence $\rPhi$ explicitly (up to an
overall sign).  For a fixed orientation on $\Sigma$, let
$\{A_i,B_i\}_{i=1}^g$ be a symplectic basis of homology classes; i.e.
$\#(A_i\cap B_j)=\delta_{i,j}$ (Kronecker $\delta$), while $\#(A_i\cap
A_j)=\#(B_i\cap B_j)=0$ for all $i,j$.  There is an induced map
$$I\colon \Wedge^k H^1(\Sigma;\Z) \longrightarrow \Wedge^k H^1(\Sigma;\Z),$$
which commutes with wedge product and satisfies $I(A_i^*)=-B_i^*$ and $I(B_i^*)=A_i^*$.
Also, there is a Hodge star operator
$$* \colon \Wedge^k H^1(\Sigma;\Z) \longrightarrow \Wedge^{2g-k} H^1(\Sigma;\Z)$$
belonging to a metric where $A_i$ and $B_i$ have length one.

Consider the map $$\phi_0\colon C\{i,j\}\longrightarrow C\{j,i\}$$
induced by $$\omega \otimes U^{-i}\mapsto (* I\omega) \otimes
U^{-j}.$$
It is easy to see that $\phi_0$ extends to an isomorphism of
$\Z$-modules, which commutes up to sign with the action by $U$ and
$\Wedge^* H_1(\Sigma;\Z)$. 

As such, it induces a graded $\Wedge^*
H_1(\Sigma;\Z)\otimes_\Z \Field[U]$-equivariant isomorphism from
$\rC\{i\geq 0\}$ to $\rC\{j\geq 0\}$.  It follows from
Lemma~\ref{lemma:AutHF} that this map is the only $\Wedge^*
H_1(\Sigma;\Z)\otimes \Field[U]$-equivariant isomorphism from
$\rC\{i\geq 0\}$ to $\rC\{j\geq 0\}$.  \qed \vskip.2cm

It will be useful to have the following terminology in place before
stating our results.

\begin{defn}
  Let $M$ be a module over $\Field[U]$, and $m\in M$. The {\em length} of $m$,
  $L_M(m)$,
  is the largest integer $\ell$ with the property that $U^\ell\cm m\neq 0$.
\end{defn}

Note that if $m\in \rC\{i,j\}\subset \Ap_s$, then
\begin{equation}
\label{eq:LengthFormula}
L_{\Ap_s}(m)=\max(i,j-s).
\end{equation}

Let $X(g,d)$ be the module over $\Z[U]$ given by
$$X(g,d)=\bigoplus_{i=0}^d \Wedge^{2g-i}H^1(\Sigma_g)\otimes_\Z
(\Z[U]/U^{d-i+1}).$$ In fact, this is a module over the ring
$\Z[U]\otimes_\Z\Wedge^* H_1(\Sigma_g)$, where the action of
$\gamma\in H_1(\Sigma_g)$ is given by
Equation~\eqref{eq:HOneAction}. We can endow $X(g,d)$ with a grading
which is lowered by one by $D_\gamma$, lowered by two by
multiplication with $U$, and centered at zero, in the sense that the
summand $\Wedge^{2g-d}H^1(\Sigma_g)$ is supported in degree zero.
Note that if $\Sym^d(\Sigma_g)$ denotes the $d$-fold symmetric product
of the surface $\Sigma_g$, then $$H_*(\Sym^d(\Sigma_g))\cong
X_{*-2d}(g,d),$$ cf.~\cite{MacDonald}.

In order to apply Proposition~\ref{prop:HomEq}, we will use coefficients
in $\Field=\Zmod{2}$.
In the interest of clarity, we restrict attention here to the case
where the Euler number $n=1$. The general case is addressed in
Theorem~\ref{thm:CircleBundlesGen} below. 

\begin{theorem}
\label{thm:CircleBundlesOne}
Let $Y$ be the circle bundle over a Riemann surface of genus $g$ with
Euler number $1$. Let $\Field=\Zmod{2}$, and consider Heegaard Floer
homology with coefficents in $\Field$.  Then, there is an isomorphism
of graded $\Field[U]$-modules
$$\HFp_{\red,*}(Y;\Field) \cong \bigoplus_{1\leq s \leq g-1} \Big(X_{*-s^2+1}(g,g-1-s) \oplus X_{*-s^2+1}(g,g-1-s)\Big).$$
\end{theorem}

\begin{proof}
  Note that $\SpinC(Y)\cong \Z^{2g}$. However, all of the homology of
  $Y$ is represented by embedded tori, and hence the adjunction
  inequality (cf. Theorem~\ref{HolDiskTwo:thm:Adjunction}
  of~\cite{HolDiskTwo}) ensures that $\HFp(Y,\spinct)$ is non-trivial
  only for the $\SpinC$ structure $\spinct$ with trivial first Chern
  class. Now, we apply Theorem~\ref{thm:Generalize} to $n=1$ surgery
  on the Borromean knot in $\#^{2g}(S^2\times S^1)$, to see that
  $\HFp(Y,\spinct)$ is isomorphic to the homology of the mapping cone
  $\Xp_{\spinct}(1)$. This is particularly simple to do, since all the
  internal differentials in $\Ap_s$ and $\Bp_s$ can be taken to
  vanish, and hence, we are left with calculating the kernel and
  cokernels of $\Dp_1$.  To simplify notation, we write simply $\Xp$
  for $\Xp_{\spinct}(1)$. In fact, we consider the reductions of all
  of these complexes modulo $2$, which we suppress from the notation;
  for example, writing simply $\Ap_s$ when we mean that this should
  be taken modulo $2$, i.e. $\Ap_s\otimes\Zmod{2}$.
  
  Observe that for all $s\geq 0$, the map $\vertp_{s}$ resp.
  $\horp_{-s}$ is surjective; and we denote its kernel by $K_s$ resp.
  $K_{-s}$.  Indeed, if $|s|$ is sufficiently large, then $K_s=0$.

  There is a chain map $\Pi_{0}\colon \Xp\longrightarrow \Ap_0$. We
  claim that this defines a surjection on homology. Specifically, any
  element $a_0\in \Ap_0$ can be extended inductively to a
  sequence $\{a_s\}_{s\in\Z}=a\in\Ker\Dp_1$ as follows: if $s\geq 0$,
  choose $a_{s+1}$ so that $\vertp(a_{s+1})=-\horp(a_s)$, and
  $a_{-s-1}$ so that $\horp(a_{-s-1})=-\vertp(a_{-s})$.  This can be
  done since $\vertp_s$ and $\horp_{-s}$ are surjective if $s\geq 0$.
  Note also that the sequence has finite support.
  
  Indeed, we claim that the projection $\Pi_0$ induces an isomorphism
  between the image of the $U^d$ on $H_*(\Xp)$ for all sufficiently
  large $d$, and $\Ap_0$, \begin{equation}
  \label{eq:HFinfIm} \left(U^d H_*(\Xp)\subset H_*(\Xp)\right)
  \stackrel{\cong}{\longrightarrow} \Ap_0. \end{equation} This can
  be seen as follows. Observe that  in all sufficiently large degrees
  (i.e. with respect to the integral grading on $\Xp$), we have
  that  for all $s>0$, both \begin{eqnarray*} \vertp_*\colon \Ap_s
  \longrightarrow \Bp_s &{\text{and}}& \horp_*\colon
  \Ap_{-s}\longrightarrow \Bp_{-s+1} \end{eqnarray*} are
  injective. It follows at once that all homogeneous elements of
  $H_*(\Xp)$ in sufficiently large degrees have non-trivial
  component in $\Ap_0$. 
  Suppose now that $a$ is such a homology class, and $a_s$ denotes its component
  in $\Ap_s$, so that $a_0$ is supported in $\rC\{i,j\}$.
  Then it is a straightforward consequence of 
  Proposition~\ref{prop:HomEq} and Equation~\eqref{eq:LengthFormula}
  that if $L_{\Ap_0}(a_0)=\max(i,j)$, then for any positive $k\geq 0$,
  $$
        L_{\Ap_s}(a_s) = 
\left\{
\begin{array}{ll}
\max(i-\frac{s(s-1)}{2},j-\frac{s(s+1)}{2}) & {\text{if $s$ is even}} \\
\max(j-\frac{s(s-1)}{2},i-\frac{s(s+1)}{2}) & {\text{if $s$ is odd}} \\
\end{array}\right.
  $$
  In particular, $L(a_s)\leq L(a_0)$ for all $s$ and hence
  $L(a_0)=L(a)$.
  It follows at once that projection $\Pi_0$ induces an injective map
  from $U^d H_*(\Xp)\subset H_*(\Xp)\longrightarrow \Ap_0$. It
  is straightforward now to see that the map $U\colon \Ap_0
  \longrightarrow \Ap_0$ is surjective. This then verifies the
  isomorphism from Equation~\eqref{eq:HFinfIm}. 

  Thus, it follows that
  for all sufficiently large $d$,  \begin{equation}
  \label{eq:ReducedHomology} \frac{H_*(\Xp)}{U^d\cm
  H_*(\Xp)}\cong \Ker \Pi_0|_{H_*(\Xp)}. \end{equation}  

  Now, for integers $i\geq 0$, let $\Xp_{|\cdot|\leq i}$ denote the
  quotient complex of $\Xp$ generated by $\Ap_{s}$ with $-i\leq s \leq
  i$ and $\Bp_{s}$ with $-i+1\leq s \leq i$.  We have projection maps
  $$ \Pi_{|\cdot|\leq i}\colon \Xp\longrightarrow \Xp_{|\cdot|\leq i},$$ 
  which, when $i=0$, coincides with the projection $\Pi_0$ considered earlier.
  Letting
  $F_i=\Ker \Pi_{|\cdot|\leq i}$, we have $$ F_0\supseteq F_1 \supseteq
  F_2 \supseteq F_3 \supseteq ... \supseteq F_i\supseteq F_{i+1}
  ... $$ Moreover, we claim that for all $s\geq 1$, there is a short
  exact sequence \begin{equation} \label{eq:ShortExactSeq} \begin{CD}
  0@>>> H_*(F_{s}) @>>> H_*(F_{s-1}) @>{\pi}>> K_s\oplus K_{-s} @>>> 0
  \end{CD} \end{equation} To see this, note that the restriction
  of $\Pi_{|\cdot|\leq s}$ to $F_{s-1}$ naturally induces a short exact sequence 
  $$ \begin{CD} 0@>>> F_{s} @>>> F_{s-1} @>>>
  \MCone\left(\vertp_s\oplus \horp_{-s}\colon \Ap_s\oplus \Ap_{-s}\longrightarrow
        \Bp_s\oplus\Bp_{-s+1}\right) @>>> 0,  \end{CD} $$ 
  and it is easy to see that the homology of the mapping cone of 
  $\vertp_s\oplus \horp_{-s}$ is identified with $K_s\oplus K_{-s}$. It remains
  to show that the map from $H_*(F_{s-1})\longrightarrow K_s\oplus K_{-s}$
  is surjective. But this follows at
  once from the fact that for all $t>s$, $\vertp_t\colon \Ap_t
  \longrightarrow \Bp_t$ and $\horp_{-t}\colon \Ap_{-t}\longrightarrow
  \Bp_{-t+1}$ are surjective. For instance, given $a_s\in K_s$ with
  $s>0$, we can use the surjectivity of $\vertp_k$ for all $k>0$ to
  inductively define a sequence $\{a_k\in \Ap_k\}_{k\in\Z}$ by
  \begin{equation} \label{eq:ExtendElement} a_k= \left\{
  \begin{array}{ll} 0 & {\text{if $k<s$}} \\ a_s & {\text{if $k=s$}}
  \\ -R_{k}\circ \horp_{k-1}(a_{k-1}) & 
   {\text{if $k>s$,}} \end{array} \right.
  \end{equation} 
  where $R_k\colon \Bp_k \longrightarrow \Ap_k$
  is a right inverse for the map $\vertp_k\colon \Ap_k\longrightarrow \Bp_k$
  (for $k\geq 0$)
  viewed as a $\Field$-module map.
  This element determines a cycle in $F_{s-1}\subset \Xp$
  whose restriction to $A_s$ is the given element $a_s\in K_s$. An
  element $a_{-s}\in K_{-s}$ can be similarly extended to a sequence
  $\{a_k\}_{k\in\Z}$ supported in $a_k$ for $k\leq -s$, by switching
  the roles of $\vertp$ and $\horp$.  From the above description of
  $C$ for the Borromean knot, it follows that $K_s\cong
  X(g,g-1-s)$. Indeed, with respect to the induced grading from $\Xp$,
  it is easy to see that $K_s\cong X_{*-s^2+1}(g,g-1-s)$. Combining
  this with the exact sequence of Equation~\eqref{eq:ShortExactSeq},
  we have proved the isomorphism claimed in the theorem as
  $\Field$-modules.

  In effect, the sequence of elements $\{a_k\}_{k\in\Z}$ defined above
  gives a right inverse to the map $\pi$ in
  Equation~\eqref{eq:ShortExactSeq}.  To verify the theorem on the
  level of $\Field[U]$-modules, it suffices to show that a careful choice
  of elements $\{a_k\}_{k\in\Z}$ provides a $\Z[U]$-equivariant
  splitting of $\pi$. Since $K_s= \rC\{i<0, j\geq s\}$, in
  view of Proposition~\ref{prop:HomEq}, the image of $\horp(K_s)$ is
  $\rC\{i\geq 0,j < -s\}$.  Now, the inclusion map of $\rC\{i\geq
  0\}\longrightarrow \rC\{\max(i,j-s-1)\geq 0\}$ provides a right
  inverse $R_s$ to $\vertp_s$ which, of course, is not a
  $\Z[U]$-module map. However, it is easy to see that its restriction
  to $\rC\{i\geq 0, j < -s\}$ is $\Z[U]$-equivariant. Moreover, this
  image is contained in the kernel of $\horp_{k+1}$. Thus, the sequence of
  elements $\{a_k\}_{k\in\Z}$ constructed from $a_s$ using this right
  inverse has at most two non-zero elements, and it is clearly a
  $\Z[U]$-equivariant splitting. When $s<0$, we use an analogous
  right inverse for $\horp_{s}$, to complete the $\Field[U]$-module
  splitting of short exact sequence 
  Equation~\eqref{eq:ShortExactSeq}. This
  completes the proof of the theorem.
\end{proof}

It is suggestive to compare the above result with Seiberg-Witten
theory over circle bundles, compare~\cite{MOY} and~\cite{Nicolaescu};
see also~\cite{KMbook} for a construction of the Floer theory in this
context.  Results from ~\cite{MOY} and~\cite{Nicolaescu} show that the
moduli space of irreducible solutions to the Seiberg-Witten equations
are identified with a disjoint union of symmetric products of the
underlying Riemann surface,
$$\coprod _{1\leq s \leq g-1} \Big(\Sym^{g-1-s}(\Sigma))\coprod
\Sym^{g-1-s}(\Sigma))\Big).$$
Theorem~\ref{thm:CircleBundlesOne} can be thought of as an algebraic
reflection of this geometric phenomenon.

According to the following result, increasing the Euler number $n$ of
the circle bundle, causes the total rank of $\HFpRed$ to drop. This
corresponds, in gauge theory, to the existence of flow-lines which
connect the vortex moduli spaces with the Jacobian torus,
compare~\cite{CircleBundles}.

Let $Y(g,n)$ denote the Euler number $n$ circle bundle over an
oriented two-manifold of genus $g$.
In this case, $\SpinC(Y(g,n))\cong \Z^{2g}\oplus \Zmod{n}$. 
By the adjunction inequality, however, only those $\SpinC$ structures
whose first Chern class is torsion have non-trivial $\HFp$.
Indeed, realizing $Y(g,n)$ as $+n$-surgery on the Borromean knot induces
an identification of $\Zmod{n}$ with these $\SpinC$ structures.

We have the following generalization of
Theorem~\ref{thm:CircleBundlesOne} which allows us to calculate the
Heegaard Floer homology for any non-trivial circle bundle over an
oriented two-manifold. Note that the hypothesis is that the Euler
number $n>0$; but the general case follows at once since $-Y(g,n)\cong
Y(g,-n)$, and hence, by the behaviour of Heegaard Floer homology under
orientation reversal, cf.~\cite{HolDisk}, we see that
$$\HFp_{\red,k}(Y(g,-n),i)\cong \HFp_{\red,-k}(Y(g,-n),i).$$

\begin{theorem}
\label{thm:CircleBundlesGen}
Let $Y(g,n)$ denote the circle bundle with Euler number $n>0$ over a
surface of genus $g$.  Then, for any choice of $[i]\in\Zmod{n}$, let
$j$ be an integer with minimal absolute value among all integers
congruent to $i\pmod{n}$.
We have that 
$$\HFp_{\red,*}(Y(g,n),i;\Zmod{2}) \cong
        \bigoplus_{\{s\equiv i\pmod{n} \big| s\neq j\}}
        X_{c(i,s)}(g,g-1-|s|),
$$
where here
$$c(i,s)=
\left\{
\begin{array}{ll}
d(n,i)-1-s+  \sum_{\{0\leq  t\leq s\big| t\equiv i\pmod{n} \}} 2t & {\text{if $s\geq 0 $}} \\
d(n,i)-1+s-  \sum_{\{s\leq t\leq 0 \big| t\equiv i\pmod{n} \}} 2t & {\text{if $s\leq 0$}} 
\end{array}
\right.
$$
\end{theorem}

\begin{proof}
This is a straightforward adaptation of the proof of
Theorem~\ref{thm:CircleBundlesOne} above.
\end{proof}

\commentable{ 
\bibliographystyle{plain} 
\bibliography{biblio} }
\end{document}